\documentclass[11pt,a4paper]{article} %

\usepackage{ifthen}
\newboolean{JMLR}
\setboolean{JMLR}{false} %

\ifthenelse{\boolean{JMLR}}{
	\usepackage{jmlr2e}
}{}

\ifthenelse{\boolean{JMLR}}{
	\usepackage[margin=0.8in, top=1in, bottom=1in]{geometry}
}{
	\usepackage[margin = 0.7in]{geometry}
}

\usepackage{amsfonts,graphicx,amsmath,amssymb,bbm,hyperref,nicefrac,enumerate,comment,mathrsfs,mathtools,enumitem,xparse,csvsimple,float,algorithm,algpseudocode,svg,subcaption,xparse,adjustbox}
\ifthenelse{\boolean{JMLR}}{
}{
	\usepackage{amsthm}
}
\definecolor{darkblue}{rgb}{0,0,0.75}
\definecolor{darkgreen}{rgb}{0,0.75,0}
\hypersetup{colorlinks=true,linkcolor=darkblue,citecolor=darkgreen}

\usepackage[titletoc,title]{appendix}
\usepackage{url}
\usepackage{booktabs}
\usepackage{rotating}

\usepackage{xcolor}
\usepackage[textwidth=0.7in]{todonotes}
\newcommand{\todoc}[1]{}

\ifthenelse{\boolean{JMLR}}{
}{
	\usepackage{cite}
}
\usepackage[sort,capitalise]{cleveref} %

\crefname{enumi}{item}{items}
\crefname{figure}{Figure}{Figures}
\crefname{equation}{}{}
\Crefname{equation}{}{}

\usepackage{tcolorbox}
\tcbuselibrary{skins, breakable, listings}  %
\newcounter{algorithmcounter}  %
\renewcommand{\thealgorithmcounter}{\arabic{algorithmcounter}}  %
\newtcolorbox[use counter=algorithmcounter]{myalgorithm}[3][]{
    enhanced,
    breakable,
    fonttitle=\bfseries,
    title=Algorithm~\thealgorithmcounter: #2,
    label={#3},  %
    label type=algorithmcounter,
    #1,
    colframe=black,
    colback=white,
    coltitle=black,
    colbacktitle=white,
    sharp corners,
    boxrule=0.5pt,
    boxsep=1mm,
    top=1mm,
    bottom=1mm,
    left=0mm,
    right=0mm
}
\crefname{algorithmcounter}{Algorithm}{Algorithms}
\Crefname{algorithmcounter}{Algorithm}{Algorithms}
\crefname{line}{line}{lines}
\Crefname{line}{Line}{Lines}

\usepackage{makecell}

\ifthenelse{\boolean{JMLR}}{}{
	\theoremstyle{plain}
	\newtheorem{theorem}{Theorem} [section] %
	
	\newtheorem{prop}[theorem]{Proposition}

	\theoremstyle{remark}
	
	\theoremstyle{definition}
	
}

\newcommand{\R}{\mathbb{R}}
\newcommand{\N}{\mathbb{N}}
\newcommand{\E}{\mathbb{E}}
\renewcommand{\P}{\mathbb{P}}

\newcommand{\cF}{\mathcal{F}}

\newcommand{\fg}{\mathfrak{g}}
\newcommand{\fx}{\mathfrak{x}}

\newcommand{\Borel}{\mathcal{B}}

\DeclarePairedDelimiterXPP{\scp}[2]{}{\langle}{\rangle}{}{#1,#2}

\DeclarePairedDelimiterXPP{\Exp}[1]{\E}{[}{]}{}{#1}

\newcommand{\norm}[1]{ \left\| #1 \right\| }

\newcommand{\normmm}[1]{{\left\vert\kern-0.25ex\left\vert\kern-0.25ex\left\vert #1 
    \right\vert\kern-0.25ex\right\vert\kern-0.25ex\right\vert}} %

\newcommand{\qand}{\qquad\text{and}}
\newcommand{\qandq}{\qquad\text{and}\qquad}
\newcommand{\andq}{\text{and}\qquad}

\newcommand{\id}{\operatorname{id}}

\makeatletter
\newcommand{\vast}{\bBigg@{3.5}}
\newcommand{\Vast}{\bBigg@{4}}
\makeatother

\DeclarePairedDelimiter{\pr}{(}{)}
\DeclarePairedDelimiter{\br}{[}{]}

\newcommand{\Ltwo}[1]{L^2(#1 ; \R)}

\ExplSyntaxOn

\NewDocumentCommand{\enum}{ O{;} m o }
 {
  \my_enum:nnn { #1 } { #2 } { #3 }
 }

\seq_new:N \l__my_enum_seq
\tl_new:N \l__my_enum_item_tl
\prop_new:N \l__verbs
\prop_put:Nnn \l__verbs {show} {shows}
\prop_put:Nnn \l__verbs {imply} {implies}
\prop_put:Nnn \l__verbs {demonstrate} {demonstrates}
\prop_put:Nnn \l__verbs {prove} {proves}
\prop_put:Nnn \l__verbs {establish} {establishes}
\prop_put:Nnn \l__verbs {ensure} {ensures}
\prop_put:Nnn \l__verbs {assure} {assures}
\int_new:N \l__number_of_args

\cs_new_protected:Nn \my_enum:nnn
 {
  \seq_set_split:Nnn \l__my_enum_seq { #1 } { #2 }
  \seq_remove_all:Nn \l__my_enum_seq {}
  \int_set_eq:NN \l__number_of_args { \seq_count:N \l__my_enum_seq }
  \seq_use:Nnnn \l__my_enum_seq { ~and~ } { ,~ } { ,~and~ }
  \IfNoValueTF{#3}{}{
    \space
    \int_compare:nNnTF{ \l__number_of_args } < {2}{ \prop_item:Nn \l__verbs {#3} }{ #3 }
  }
 }

\ExplSyntaxOff

\DeclareFontEncoding{LS1}{}{}
\DeclareFontSubstitution{LS1}{stix}{m}{n}
\DeclareMathAlphabet{\mathscr}{LS1}{stixscr}{m}{n}

\newcommand{\multdim}[2]{{#1}}

\newcommand{\providecommandordefault}[2]{%
    \providecommand{#1}{}%
    \renewcommand{#1}{#2}%
}

\makeatletter
\ExplSyntaxOn
\seq_new:N \g_abbrs
\prop_new:N \g_abbr_counts
\tl_new:N \l_abbr_count_tl

\NewDocumentCommand{\abbr}{m m O{#1} m m O{#4}}{
    \expandafter\newcommand\csname#3\endcsname[1][]{
        \seq_if_in:NnTF \g_abbrs {#1} {
            \prop_get:NnN \g_abbr_counts {#1} \l_abbr_count_tl
            \prop_gput:Nnx \g_abbr_counts {#1} {\int_eval:n {\l_abbr_count_tl + 1}}
            \hyperref[#1]{#1}
        } {
            \seq_gput_left:Nn \g_abbrs {#1}
            \prop_gput:Nnn \g_abbr_counts {#1} {1}
            \expandafter\gdef\csname#1@def\endcsname{#2}
            \phantomsection\label{#1}
            \str_if_eq:nnTF{##1}{}{\emph{#2}}{##1}~(\hyperref[#1]{#1})
        }
    }
    \expandafter\newcommand\csname#6\endcsname[1][]{
        \seq_if_in:NnTF \g_abbrs {#1} {
            \prop_get:NnN \g_abbr_counts {#1} \l_abbr_count_tl
            \prop_gput:Nnx \g_abbr_counts {#1} {\int_eval:n {\l_abbr_count_tl + 1}}
            \hyperref[#1]{#4}
        } {
            \expandafter\gdef\csname#1@def\endcsname{#5}
            \seq_gput_left:Nn \g_abbrs {#1}
            \prop_gput:Nnn \g_abbr_counts {#1} {1}
            \phantomsection\label{#1}
            \str_if_eq:nnTF{##1}{}{\emph{#5}}{##1}~(\hyperref[#1]{#4})
        }
    }
}

\ExplSyntaxOff
\makeatother

\abbr{ANN}{artificial neural network}{ANNs}{artificial neural networks}
\abbr{SGD}{Stochastic Gradient Descent}{SGDs}{Stochastic Gradient Descents}
\abbr{PDE}{partial differential equation}{PDEs}{partial differential equations}
\abbr{ODE}{ordinary differential equation}{ODEs}{ordinary differential equations}
\abbr{FNO}{Fourier neural operator}{FNOs}{Fourier neural operators}
\abbr{ADANN}{Algorithmically Designed Artificial Neural Networks}{ADANNs}{Algorithmically Designed Artificial Neural Networks}
\abbr{LIRK}{linearly implicit Runge-Kutta}{LIRKs}{Linearly Implicit Runge-Kutta}
\abbr{GELU}{Gaussian Error Linear Unit}{GELUs}{Gaussian Error Linear Units}
\abbr{DeepONet}{deep operator network}{DeepONets}{deep operator networks}

\newcommand{\nrEvalRuns}

\providecommandordefault{\nrspacediscConcrete}{64}
\providecommandordefault{\nrValSamplesConcrete}{2^{14}}
\providecommandordefault{\valNrTimesteps}{1500}
\providecommandordefault{\valNrSpacesteps}{512}

\providecommandordefault{\nrTrainSamplesConcrete}{2^{18}}
\providecommandordefault{\trainNrTimesteps}{1000}
\providecommandordefault{\trainNrSpacesteps}{256}

\providecommandordefault{\lrTestTrainsteps}{50}

\providecommandordefault{\nrEvalRuns}{1000}

\providecommandordefault{\nrRunsOL}{5}

\newcommand{\SGonedT}{2}
\newcommand{\SGonedSpaceSize}{1}
\newcommand{\SGonedLaplaceFactor}{\frac{1}{100}}

\newcommand{\SGoneddecayRate}{4}
\newcommand{\SGonedoffset}{10^{2.5}}
\newcommand{\SGonedvar}{10^{10}}

\newcommand{\SGonedInitDistr}{\mathcal{N}(0, \SGonedvar(\SGonedoffset \id_{\initialValues} - \Delta_x)^{-\SGoneddecayRate})}

\newcommand{\SGonedSpaceStep}{64}
\newcommand{\SGonedParams}{(0.1, 1.2) \times (0.25, 1.2)}
\newcommand{\SGonedBatchSize}{256}
\newcommand{\SGonedILRsteps}{50}
\newcommand{\SGonedEvalSteps}{400}
\newcommand{\SGonedImprTol}{0.96}
\newcommand{\SGonedNrOLRuns}{12}
\newcommand{\SGonedNrGridRuns}{25}
\newcommand{\SGonedNrOptRuns}{12}
\newcommand{\SGonedRefAlgOne}{Spectral/Crank-Nicolson}
\newcommand{\SGonedRefAlgTwo}{explicit midpoint}
\newcommand{\SGonedNrTrainSamples}{2^{18}}
\newcommand{\SGonedNrValidSamples}{2^{14}}
\newcommand{\SGonedNrTestSamples}{2^{14}}
\newcommand{\SGonedNrTrainSpaceDiscr}{256}
\newcommand{\SGonedNrValidSpaceDiscr}{512}
\newcommand{\SGonedNrTestSpaceDiscr}{512}
\newcommand{\SGonedNrTrainTimeSteps}{1000}
\newcommand{\SGonedNrValidTimeSteps}{1500}
\newcommand{\SGonedNrTestTimeSteps}{1500}

\newcommand{\SGtwodSpaceStep}{32}
\newcommand{\SGtwodParams}{\SGonedParams}
\newcommand{\SGtwodBatchSize}{128}
\newcommand{\SGtwodILRsteps}{50}
\newcommand{\SGtwodEvalSteps}{400}
\newcommand{\SGtwodImprTol}{0.96}
\newcommand{\SGtwodNrOLRuns}{12}
\newcommand{\SGtwodNrGridRuns}{-}
\newcommand{\SGtwodNrOptRuns}{12}
\newcommand{\SGtwodRefAlgOne}{Spectral/Crank-Nicolson}
\newcommand{\SGtwodRefAlgTwo}{explicit midpoint}
\newcommand{\SGtwodNrTrainSamples}{2^{16}}
\newcommand{\SGtwodNrValidSamples}{2^{11}}
\newcommand{\SGtwodNrTestSamples}{2^{11}}
\newcommand{\SGtwodNrTrainSpaceDiscr}{64}
\newcommand{\SGtwodNrValidSpaceDiscr}{128}
\newcommand{\SGtwodNrTestSpaceDiscr}{128}
\newcommand{\SGtwodNrTrainTimeSteps}{1000}
\newcommand{\SGtwodNrValidTimeSteps}{1500}
\newcommand{\SGtwodNrTestTimeSteps}{1500}

\newcommand{\BurgersT}{1}
\newcommand{\BurgersSpaceSize}{2\pi}
\newcommand{\BurgersLaplaceFactor}{\frac{1}{10}}

\newcommand{\BurgersdecayRate}{6}
\newcommand{\Burgersoffset}{10}
\newcommand{\Burgersvar}{10^6}

\newcommand{\BurgersInitDistr}{\mathcal{N}(0, \Burgersvar(\Burgersoffset \id_{\initialValues} - \Delta_x)^{-\BurgersdecayRate})}

\newcommand{\BurgersSpaceStep}{32}
\newcommand{\BurgersParams}{(0.1, 1)^2}
\newcommand{\BurgersBatchSize}{1024}
\newcommand{\BurgersILRsteps}{50}
\newcommand{\BurgersEvalSteps}{400}
\newcommand{\BurgersImprTol}{0.96}
\newcommand{\BurgersNrOLRuns}{16}
\newcommand{\BurgersNrGridRuns}{16}
\newcommand{\BurgersNrOptRuns}{-}
\newcommand{\BurgersRefAlgOne}{Spectral/Crank-Nicolson}
\newcommand{\BurgersRefAlgTwo}{explicit midpoint}
\newcommand{\BurgersNrTrainSamples}{2^{18}}
\newcommand{\BurgersNrValidSamples}{2^{14}}
\newcommand{\BurgersNrTestSamples}{2^{14}}
\newcommand{\BurgersNrTrainSpaceDiscr}{128}
\newcommand{\BurgersNrValidSpaceDiscr}{256}
\newcommand{\BurgersNrTestSpaceDiscr}{256}
\newcommand{\BurgersNrTrainTimeSteps}{1000}
\newcommand{\BurgersNrValidTimeSteps}{1500}
\newcommand{\BurgersNrTestTimeSteps}{1500}

\newcommand{\BurgersDMarch}{(\BurgersSpaceStep, 256, 1024, 256, \BurgersSpaceStep)}

\newcommand{\ReactDiffT}{1}
\newcommand{\ReactDiffSpaceSize}{2}
\newcommand{\ReactDiffLaplaceFactor}{\frac{5}{100}}
\newcommand{\ReactDiffReactionRate}{2}

\newcommand{\ReactDiffdecayRate}{4}
\newcommand{\ReactDiffoffset}{100}
\newcommand{\ReactDiffvar}{10^8}

\newcommand{\ReactDiffInitDistr}{\mathcal{N}(0, \ReactDiffvar(\ReactDiffoffset \id_{\initialValues} - \Delta_x)^{-\ReactDiffdecayRate} - 0.8 \id_{\initialValues})}

\newcommand{\ReactDiffSpaceStep}{128}
\newcommand{\ReactDiffParams}{(0.1, 1.3)^2}
\newcommand{\ReactDiffBatchSize}{256}
\newcommand{\ReactDiffILRsteps}{50}
\newcommand{\ReactDiffEvalSteps}{400}
\newcommand{\ReactDiffImprTol}{0.97}
\newcommand{\ReactDiffNrOLRuns}{16}
\newcommand{\ReactDiffNrGridRuns}{16}
\newcommand{\ReactDiffNrOptRuns}{-}
\newcommand{\ReactDiffRefAlgOne}{FDM/Crank-Nicolson}
\newcommand{\ReactDiffRefAlgTwo}{explicit midpoint}
\newcommand{\ReactDiffNrTrainSamples}{2^{18}}
\newcommand{\ReactDiffNrValidSamples}{2^{14}}
\newcommand{\ReactDiffNrTestSamples}{2^{14}}
\newcommand{\ReactDiffNrTrainSpaceDiscr}{512}
\newcommand{\ReactDiffNrValidSpaceDiscr}{1024}
\newcommand{\ReactDiffNrTestSpaceDiscr}{1024}
\newcommand{\ReactDiffNrTrainTimeSteps}{1000}
\newcommand{\ReactDiffNrValidTimeSteps}{1500}
\newcommand{\ReactDiffNrTestTimeSteps}{1500}

\newcommand{\ReactDiffDMarch}{(\ReactDiffSpaceStep, 512, 1024, 512, \ReactDiffSpaceStep)}

\begin{document}

\title{
	Algorithmically Designed Artificial \\ Neural Networks (ADANNs): \\
	Higher order deep operator learning for \\ parametric partial differential equations
}

\ifthenelse{\boolean{JMLR}}{
	\ShortHeadings{Algorithmically Designed Artificial Neural Networks}{Jentzen, Riekert, von Wurstemberger}
	\firstpageno{1}

	\author{\name 
		Arnulf Jentzen
		\email ajentzen@cuhk.edu.cn, ajentzen@uni-muenster.de\\
		\addr School of Data Science and Shenzhen Research Institute of Big Data \\
		The Chinese University of Hong Kong, Shenzhen (CUHK-Shenzhen) \\
		China \\
		\addr  
		Applied Mathematics: Institute for Analysis and Numerics\\
		University of M\"unster\\
		Germany
		\AND
		\name Adrian Riekert \email ariekert@uni-muenster.de \\
		\addr 
		Applied Mathematics: Institute for Analysis and Numerics\\ 
		University of M\"unster\\ 
		Germany
		\AND
		\name Philippe von Wurstemberger
		\email philippevw@cuhk.edu.cn, vwurstep@ethz.ch\\
		\addr 
		School of Data Science \\ 
		The Chinese University of Hong Kong, Shenzhen (CUHK-Shenzhen) \\
		China \\
		\addr
		Risklab, Department of Mathematics\\
		ETH Zurich\\
		Switzerland
	}

	\editor{...}
}{

	\author{
	Arnulf Jentzen$^{1,2}$,
	Adrian Riekert$^{3}$,
	and
	Philippe von Wurstemberger$^{4, 5}$
	\bigskip\\
	\small{$^1$ School of Data Science and School of Artificial Intelligence, 
	The Chinese University} \vspace{-0.1cm}\\
	\small{of Hong Kong, Shenzhen (CUHK-Shenzhen), China; e-mail: \texttt{ajentzen}\textcircled{\texttt{a}}\texttt{cuhk.edu.cn}}\smallskip\\
	\small{$^2$ Applied Mathematics: Institute for Analysis and Numerics,}\vspace{-0.1cm}\\
	\small{University of M\"unster, Germany; e-mail: \texttt{ajentzen}\textcircled{\texttt{a}}\texttt{uni-muenster.de}}\smallskip\\
	\small{$^3$ Applied Mathematics: Institute for Analysis and Numerics,}\vspace{-0.1cm}\\
	\small{University of M\"unster, Germany; e-mail: \texttt{ariekert}\textcircled{\texttt{a}}\texttt{uni-muenster.de}}\smallskip\\
	\small{$^4$ School of Data Science, The Chinese University of Hong Kong,} \vspace{-0.1cm}\\
	\small{Shenzhen (CUHK-Shenzhen), China; e-mail: \texttt{philippevw}\textcircled{\texttt{a}}\texttt{cuhk.edu.cn}}\smallskip\\
	\small{$^5$ Risklab, Department of Mathematics, ETH Zurich, Switzerland} 
	}
}

\maketitle

\begin{abstract}
 
	In this article we propose a new deep learning approach to approximate operators related to 
	parametric 
	partial differential equations (PDEs).
	In particular, we introduce a new strategy to design specific artificial neural network (ANN) architectures in conjunction with specific ANN initialization schemes 
	which are tailor-made for the particular approximation problem under consideration. 
	In the proposed approach we combine 
	efficient classical numerical approximation techniques 
	with deep operator learning methodologies.
	Specifically, 
	we introduce customized adaptions of existing ANN architectures 
	together with specialized initializations for these ANN architectures 
	so that at initialization we have that 
	the ANNs closely mimic a chosen 
	efficient classical numerical algorithm 
	for the considered approximation problem.
	The obtained ANN architectures and their initialization schemes 
	are thus strongly inspired by numerical algorithms 
	as well as by popular deep learning methodologies from the literature
	and in that sense we refer to the introduced ANNs in conjunction 
	with their tailor-made initialization schemes as 
	\emph{Algorithmically Designed Artificial Neural Networks} (ADANNs).  
	We numerically test the proposed ADANN methodology in the case 
	of 
	several
	parametric 
	PDEs.
	In the tested numerical examples
	the ADANN methodology significantly outperforms 
	existing classical approximation algorithms as well as 
	existing deep operator learning methodologies from the literature.

\end{abstract}

\ifthenelse{\boolean{JMLR}}{
}{
	\tableofcontents
}

\section{Introduction}

Deep learning approximation methods -- usually consisting of deep \ANN\ models trained 
through \SGD\ optimization methods -- are nowadays among the
most heavily employed approximation methods in the digital world. 
They are behind most of the recent success in artificial intelligence and machine learning in areas such as computer vision (cf., e.g, \cite{Krizhevsky2012,Dosovitskiy2021,Voulodimos2018} and the references therein), 
natural language processing (cf., e.g., \cite{Vaswani2017,Brown2020,OpenAI2023,Qiu2020,Rogers2020} and the references therein), 
and 
speech recognition  (cf., e.g., \cite{Hinton2012,Deng2013,Yu2015,Mehrish2023} and the references therein)
and are increasingly used in many other fields.
One such field is scientific computing where in recent years  
deep learning technologies have been intensively applied to various problems including to 
the numerical approximation of 
\PDEs.

In particular, deep learning approximation methods have been developed to 
approximately solve high-dimensional nonlinear \PDEs\ 
(see, e.g., \cite{han2018solving,Weinan2017,blechschmidt2021three,beck2020overviewPublished,Weinan2020a,Germain2021,sirignano2018dgm,Jentzen2023,Nusken21}
and the references therein)
such as 
high-dimensional nonlinear pricing problems from financial engineering and 
Hamiltonian-Jacobi-Bellman equations from optimal control. 
In the context of 
such 
high-\discretionary{dimen-}{sional}{dimensional}
nonlinear 
\PDEs\, the progress achieved by 
deep learning approximation methods 
is obvious as there are -- 
except in some special cases 
(see, e.g., \cite{HenryLabordere12,henry2019branching,Nguwi2022} and the references therein for Branching-type methods
and see, e.g., \cite{E2019,EHutzenthalerJentzenKruse16published,Weinan2020a,Hutzenthaler_2020} and the references therein for multilevel Picard methods) 
-- 
essentially no practical alternative numerical approximation methods.
The striking feature of deep learning methods in the context of such high-dimensional problems is that in many situations 
numerical simulations suggest 
that the computational effort of such methods 
only grows at most polynomially 
in the input dimension
of the problem under consideration. 
In contrast, classical numerical methods usually suffer under the so-called 
\emph{curse of dimensionality} 
(cf., e.g., \cite{Bellman66}, \cite[Chapter 1]{NovakWozniakowski08}, and \cite[Chapter 9]{Novak08})
in the sense that 
the computational effort grows at least exponentially in the dimension.

There is also a vast literature on deep learning approximation methods 
for low-dimensional \PDEs\ 
(cf., e.g., \cite{Karniadakis21,Raissi19,Jentzen2023,beck2020overviewPublished} and the references therein).
For 
low-dimensional \PDEs\, in most cases,
there usually already exist a number 
of efficient classical (non-deep learning based) approximation methods
in the scientific literature (cf., e.g., \cite{tadmor2012review,Jovanovic2014,LeVeque2007,Bartels2015}).
Nonetheless, there are several convincing arguments 
that deep learning approximation methods 
might have the potential to significantly outperform 
such traditional approximation methods 
from classical numerics. 
One situation where this strongly applies 
is in the context of \emph{parametric}
\PDE\ approximation problems. 
Specifically,  
in applications one is often not only interested 
in approximately solving a 
\PDE\ model once but 
instead there is often the need to approximately solve a \PDE\ model
repeatedly but with different 
initial values and/or different model parameters.
The idea of deep learning approaches in this context 
is to try to not only solve one fixed \PDE\ but instead 
to learn the whole solution 
operator
which maps initial values and model parameters 
to corresponding \PDE\ solutions. 
Even though the original \PDE\ model has often only 
one to three space-dimensions, 
the associated approximation problem for the solution operator 
typically becomes very high-dimensional due to the high number 
of parameters required to approximately describe the initial value 
and the model parameters.
Because of their apparent capacity to overcome the curse of dimensionality, deep learning methods, therefore
seem to be very natural candidates for such approximation problems.
Deep learning methods in this situation are then 
often referred to as \emph{(deep) operator learning} approaches (cf., e.g., \cite{Li2021,Li2020,Lu2021}). 
However, even though very remarkable 
advances have been accomplished in this area of research, 
for instance, by means of 
so-called \FNOs\
(see \cite{Li2021}), 
so far in most situations deep operator learning techniques do not outperform 
the most efficient higher order classical numerical methods.
This is also not entirely surprising due to fundamental lower bounds 
established in the literature 
that a wide class of methods, including typical deep learning approximations,
can in general not overcome the curse of dimensionality 
in the $L^\infty$-norm 
(cf., e.g., \cite{hs99,Heinrich2006,Grohs2023published}).

It is precisely the objective of this work to introduce a new operator 
learning approach 
which aims to overcome this challenge 
by combining efficient classical numerical methods
with deep operator learning techniques.
In particular, we introduce a new strategy 
to design specific \ANN\ architectures in conjunction with specific \ANN\ initialization schemes
which are tailor-made for the particular approximation problem under consideration.
The obtained \ANN\ architectures and their initialization schemes
are strongly inspired by numerical algorithms
as well as by popular deep learning methodologies from the literature
and in that sense we refer to the introduced \ANN\ architectures in conjunction
with their tailor-made initialization schemes as
\ADANNs.
We numerically test the \ADANN\ methodology proposed in this paper in the case 
of several parametric \PDEs.
In the tested numerical examples
the \ADANN\ methodology significantly outperforms classical approximation algorithms as well as existing deep operator learning methodologies from the literature.

We now briefly describe some key aspects of the \ADANN\ methodology proposed in this paper in more detail.
The architecture of \ADANNs\ has two components which are added to each other and trained separately:
a \emph{base model} and a \emph{difference model}.
A base model is designed together with highly specialized initializations for that model such that at those initializations the base model exactly represents a family of efficient classical numerical algorithms.
It is then trained to approximate the considered target operator using \SGD-type methods starting at its highly specialized initializations. 
Loosely speaking, the training of a base model can be seen as improving the classical numerical algorithms.
On the other hand, a difference model can be chosen to be any operator learning model from the literature
and is trained to approximate the difference between the target operator and a base model.
A further key component of the \ADANN\ methodology is the use an additional optimization approach over several training runs for the base and difference model starting from different initializations.

We now discuss some ideas in the scientific literature which are related to the \ADANN\ methodology introduced in this paper. 
We roughly divide the related literature into three categories.
First, we consider approaches similar to the idea of base models in the \ADANN\ methodology, in the sense that existing numerical algorithms are somehow made trainable and are subsequently improved by means of \SGD-type methods.
One such approach -- by which the \ADANN\ methodology was partially inspired -- is the \emph{learning the random variable} methodology in \cite{Becker2024}
where \emph{Monte Carlo neural networks} have been introduced, which have the property that at suitable initializations the realizations of those networks correspond to sample realizations of Monte Carlo algorithms.
Several approaches where Runge-Kutta methods for \ODEs\ are improved by considering parts of the Butcher tableau as trainable parameters can, e.g., be found in \cite{Tsitouras2002,Anastassi2014,Dehghanpour2015,Ouala2021,Mishra2019}.
Moreover, in \cite{Mishra2019} also certain parameters within numerical algorithms for \PDEs\ are considered as trainable parameters.
While the approaches in \cite{Tsitouras2002,Anastassi2014,Dehghanpour2015,Ouala2021,Mishra2019} are similar to the idea of base models in the \ADANN\ methodology (and can be considered special cases of the \ADANN\ methodology), they usually only train a handful of parameters which have an inherent meaning within the considered algorithms and in most cases even enforce some type of order conditions on the trained parameters, whereas we use a classical numerical algorithm as a starting point to design a base model but then do not restrict the model to stay within a class of known algorithms.
In \cite{Katrutsa2017,Greenfeld2019} the restriction and prolongation operators in multigrid methods are considered as trainable parameters.
In \cite{Wang1998} \ANNs\ are designed to emulate Runge-Kutta methods, not to learn solution operators, but to learn the dynamics of unknown \ODEs\ from observed trajectories.
In \cite{BarSinai2019} finite difference approximations of spatial derivative operators for solutions of spatio-temporal \PDEs\ are improved by means of \SGD-type methods and subsequently used to perform time-integrations of the considered \PDEs.
Similarly, \cite{Kossaczka2023} improve finite difference approximations of spatial derivatives in spatio-temporal \PDEs\, but instead of directly improving the finite difference coefficients like in \cite{BarSinai2019} they learn the truncation error by means of \ANNs.

Second, we mention approaches which are similar to the idea of difference models in the \ADANN\ methodology, in the sense that \ANNs\ are used to improve classical numerical algorithms by learning the residual between a classical algorithm and the corresponding target quantity.
There are a number of approaches in the literature where \ANNs\ are trained to approximate the difference between coarse-grid and fine-grid approximations of \PDE\ solutions, notably within the context of large eddy simulations for the simulation of turbulent flows, cf., e.g., \cite{Dresdner2022,Kochkov2021,List2022,Subel2021,Maulik2019,San2018,Frezat2022,Huang2023a} and the references therein.
Moreover, there are also approaches where \ANNs\ are used to learn the time stepping error of \ODE\ integration schemes, cf., e.g., \cite{Shen2020,Huang2023,Song2024}.

Third, we list some approaches which aim to combine classical methods with deep learning techniques but are neither close to the idea of base models nor to the idea of difference models in the \ADANN\ methodology.
A number of approaches use \ANNs\ to guide classical methods by training the \ANNs\ to identify regions where the classical methods need to be adjusted, cf., e.g., \cite{Ray2018,Discacciati2020,Fidkowski2021,Bohn2021}.
In \cite{Hsieh2019} traditional iterative solvers are improved by adding some learnable parameters to the iterator.
In \cite{Tompson2017} a convolutional neural network is used to efficiently approximate a linear projection which is then used as a component in a standard solver for fluid flows.
In \cite{Brevis2021} \ANNs\ are used to approximate optimal test spaces in the context of finite element methods.

Next we mention several other promising deep operator learning approaches in the literature, all of which can, in principle, be used as difference models within the \ADANN\ methodology.
One of the most successful methods in practice are the \FNOs\ introduced in \cite{Li2021}.
The derivation of \FNOs\ is based on \cite{Li2020}, an earlier paper by the same authors, where so-called graph kernel networks are employed. 
In \cite{Li2022} \FNOs\ are generalized to more complicated geometries.
In \cite{Brandstetter2022} the \FNO\ methodology is extended by using Clifford layers, where calculations are performed in higher-dimensional non-commutative Clifford algebras.
Another successful approach is the \DeepONet\  architecture introduced in \cite{Lu2021}, which consists of two types of \ANNs\ that take as input the output space points and the input function values, respectively.
For a comparison between the \DeepONet\ and \FNO\ methodologies we refer to \cite{Lu2022}.
In \cite{Lanthaler2022} \DeepONets\ are generalized to a more sophisticated nonlinear architecture.
In \cite{Pham2022} operators on Wasserstein spaces, for example, mean-field interactions of measures, are learned using networks based on standard \ANNs\ and \DeepONets.
In \cite{Nelsen2021} operators between Banach spaces are approximated by using random feature maps associated to operator-valued kernels.
In \cite{Liu2022} the entire flow map associated to \ODEs\ is approximated by training a different \ANN\ in each time-step and combining these \ANNs\ with classical Runge-Kutta methods on different time scales.
For approaches to build operator learning architectures based on convolutional neural networks, we refer, e.g., to \cite{Guo2016,Raonic2023,Heiss2023,Zhu2018,Khoo2021}.
We also refer to \cite{Kovachki2021,Chen2023} for estimates for approximation and generalization errors in network-based operator learning for \PDEs.
Finally, we refer to \cite[Appendix D]{Brandstetter2022} and \cite[Section 1.7.4]{Jentzen2023}
for more detailed literature overviews on operator learning approaches.

The remainder of this article is organized as follows. 
In \cref{sect:rough_overview} we introduce the main ideas of the \ADANN\ methodology in an abstract setting.
In \cref{sect:semilinear_heat} we describe in detail a specific base model design in the case of semilinear heat \PDEs.
Finally, in \cref{sect:simul} we present four numerical simulations comparing the \ADANN\ methodology to classical methods and operator learning methods from the literature.

\section{Overview of the ADANN methodology}
\label{sect:rough_overview}

\newcommand{\adann}{\mathscr{A}}

\newcommand{\solOp}{\mathcal{S}}
\newcommand{\solOpAlt}{\tilde{\mathcal{S}}}
\newcommand{\initialValues}{\mathcal{I}}
\newcommand{\EndValues}{\mathcal{O}}
\newcommand{\propAlg}{\Phi}
\newcommand{\refAlg}{\Psi}
\newcommand{\nrinitdiscr}{I}
\newcommand{\nrenddiscr}{O}
\newcommand{\params}{\mathfrak{P}}
\newcommand{\nrbaseParams}{{\mathbf{d}_{\text{Base}}}}
\newcommand{\initParams}{\mathbf{W}}
\newcommand{\adannBase}{\mathscr{B}}
\newcommand{\nrtimesteps}{{M}}
\newcommand{\cost}{\delta}
\newcommand{\loss}{L_\text{Base}}
\newcommand{\intLoss}{\mathbf{L}^{(\text{Base})}}
\newcommand{\nrReftimesteps}{\mathbf{R}}
\newcommand{\initialRV}{\mathfrak{I}}
\newcommand{\IVvariable}{g}
\newcommand{\baseSGD}{\mathcal{W}}
\newcommand{\diffSGD}{\Theta}
\newcommand{\batchsize}{B}

\newcommand{\approxParams}{\mathfrak{W}}
\newcommand{\diffann}{\mathscr{D}}
\newcommand{\nrdiffParams}{{\mathbf{d}_{\text{Diff}}}}
\newcommand{\estLoss}{\varepsilon}
\newcommand{\diffLoss}{L_\text{Diff}}
\newcommand{\differenceintLoss}[1]{\mathbf{L}^{(\text{Diff})}_{#1}}
\newcommand{\EndLoss}{\mathfrak{e}}

\newcommand{\error}{\mathfrak{E}}

\newcommand{\argmin}{\text{argmin}}

\DeclarePairedDelimiterX{\lossmetric}[1]{\lVert}{\rVert}{#1}

In this section we describe the \ADANN\ methodology in an abstract setting.
For this, we consider the problem of numerically approximating a measurable operator
\begin{equation}
	\label{base_model:eq1}
	\begin{split} 
		\solOp \colon \initialValues \to \EndValues
	\end{split}
\end{equation}
where $\initialValues$ and $\EndValues$ are topological vector spaces\footnote{
	For example, we think of the mapping which assigns to all suitable initial values of a \PDE\ of evolutionary type, the corresponding terminal value of the \PDE.
}.
We assume that we are given a random input variable
$
	\initialRV \colon \Omega \to \initialValues
$
on a probability space 
$(\Omega, \cF, \P)$
and a continuous seminorm
$
	\lossmetric{\cdot} \colon \EndValues \to [0,\infty)
$
on the output space.
We measure the quality of a measurable approximation $\solOpAlt \colon \initialValues \to \EndValues$ of the operator in \eqref{base_model:eq1} by means of the $L^2$-error
$
	\error(\solOpAlt)
\in 
	[0, \infty] 
$ given by
\begin{equation}
	\label{T_B_D}
	\begin{split}
		\error(\solOpAlt)
		=
		\pr[\big]{
			\Exp[\big]{
				\lossmetric{\solOpAlt(\initialRV) - \solOp(\initialRV)}^2
			}
		}^{1/2}.
	\end{split}
\end{equation}

\begin{figure}
	\centering
	\includegraphics[width=1\linewidth]{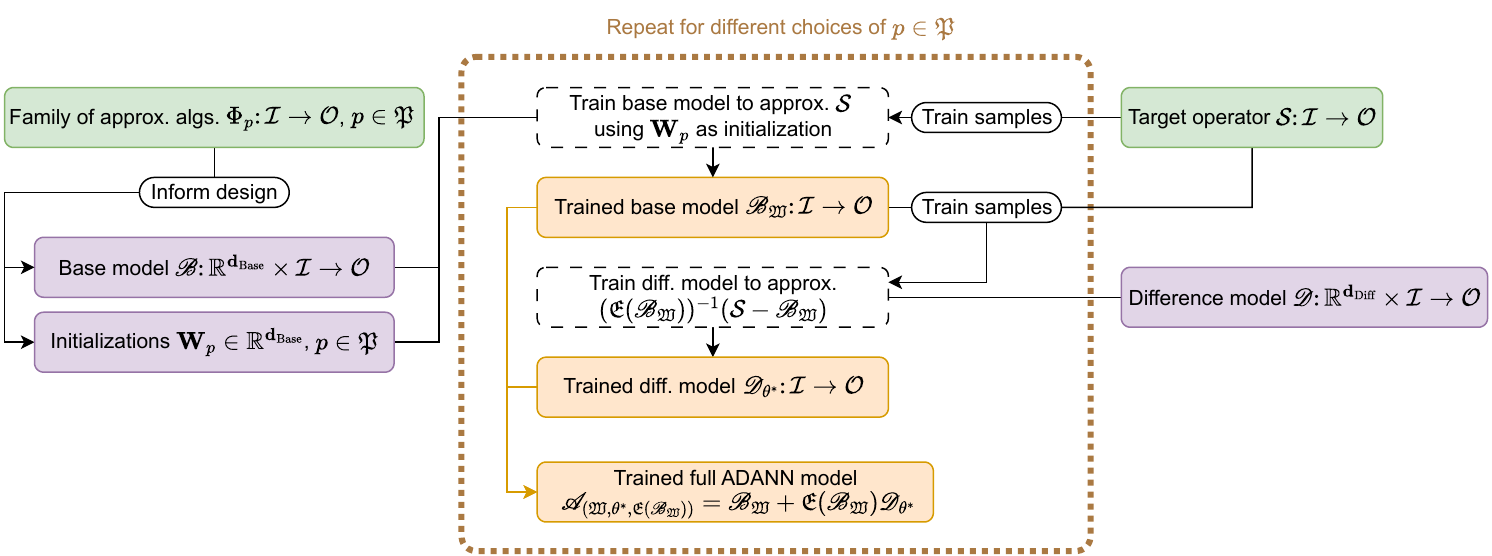}
	\caption{\label{fig:adann_schematic}Schematic overview of the \ADANN\ methodology developed in \cref{sect:rough_overview}.}
\end{figure}

Roughly speaking, the \ADANN\ methodology
proposed in this paper relies on the following three main ingredients to produce an efficient approximation of the operator $\solOp$:
\begin{enumerate}[label=(\roman *)]
	\item \label{adanns:item1}
		\emph{Base model with highly specialized initializations}:
		Based on a family of classical 
		numerical approximation algorithms we design a tailor-made problem-specific \ANN-type model together with a family of highly specialized initializations for the model
		and train the model to approximate $\solOp$ (cf.\ \cref{sect:adann_base}).
	      
	\item \label{adanns:item2}
	      \emph{Difference model}:
	      We employ an existing operator learning model from the literature to approximate the difference between realizations of the base model of step \ref{adanns:item1} and the operator $\solOp$. Adding the difference model to the base model results in the \emph{full \ADANN\ model} (cf.\ \cref{sect:difference_training}). 
	      
	\item \label{adanns:item3}
	      \emph{Optimization over base model initializations}: 
	      We repeat the training of the full \ADANN\ model with different (possibly random) highly specialized initializations of the base model and different standard random initializations of the difference model. The highly specialized initializations of the base model are chosen based on a suitable additional optimization approach aiming to minimize the error of the best full \ADANN\ model over all training runs (cf.\ \cref{sect:random_init}).
	      
\end{enumerate}
The combination of the above three components into the \ADANN\ methodology is detailed in \cref{sect:adann_pseudocode} in the form of pseudocodes and illustrated in \cref{fig:adann_schematic}.

\subsection{Base model with highly specialized initializations}
\label{sect:adann_base}

Constructing a base model requires a parametric family 
$\propAlg_p \colon \initialValues \to \EndValues$, $p \in \params$,
of classical numerical algorithms\footnote{For example, we think of a family of Runge-Kutta and/or finite element methods.}
indexed over a parameter set $\params$
which approximate the operator $\solOp$ in the sense that 
for every $p \in \params$ we have that
\begin{equation}
	\label{base_model:eq2}
	\begin{split} 
		\propAlg_p \approx \solOp.
	\end{split}
\end{equation}
A \emph{base model}
corresponding to the \emph{designing algorithms} in  \eqref{base_model:eq2}
with 
$
	\nrbaseParams \in \N = \{ 1, 2, 3, \dots \}
$
trainable parameters 
is then given by a model 
$
	\adannBase
	=
	(\adannBase_W(\IVvariable))_{(W, \IVvariable) \in \R^\nrbaseParams \times \initialValues}
	\colon
	\R^\nrbaseParams \times \initialValues \to \EndValues
$
which can reproduce all the algorithms in \eqref{base_model:eq2}. %
In mathematical terms, we require that
for every
$p \in \params$
there exist model parameters
$\initParams_p \in \R^\nrbaseParams$
with 
\begin{equation}
	\label{base_model:eq4}
	\begin{split} 
		\adannBase_{\initParams_p}
		=
		\propAlg_p
		\approx 
		\solOp.
	\end{split}
\end{equation}
Note that \eqref{base_model:eq4} implies that we already have parameters for the base model which yield reasonable approximations of $\solOp$. 
In the training of the base model we suggest to improve these approximations.
Specifically, we propose to minimize the loss
$\intLoss \colon \R^\nrbaseParams\to [0,\infty]$
satisfying\footnote{In applications, the target operator $\solOp$ is typically not known and has to be approximated for the training process.
	Such approximations can, for example, be obtained by using the same class of algorithms as the designing algorithms in \eqref{base_model:eq2} but with much finer discretizations.} for all 
$W \in \R^\nrbaseParams$
that
\begin{equation}
	\label{base_model:eq6}
	\begin{split} 
		\intLoss(W)
	=
		\Exp*{
			\lossmetric{\adannBase_W(\initialRV) - \solOp(\initialRV)}^2
		}
	=
		\error(\adannBase_W)^2
	\end{split}
\end{equation}
by means of \SGD-type processes initialized with model parameters from the set $\{\initParams_p \colon p \in \params \}$.

For concrete examples of base models we refer to \cref{sect:semilinear_heat} where we derive designing algorithms and corresponding base models in the context of semilinear heat \PDEs\ and 
to \cref{sect:SG_base,sect:Burgers_base,sect:RD_base} where we specify the base models used in our numerical simulations.

\subsection{Difference models}
\label{sect:difference_training}
\newcommand{\nrruns}{R}

To combine a base model from \cref{sect:adann_base} with deep operator learning models from the literature we propose to employ an existing deep operator learning model to approximate the difference between the operator $\solOp$ and the base model.
For this, we introduce a \emph{difference model}
$
	\diffann = (\diffann_\theta(\IVvariable))_{(\theta, \IVvariable) \in \R^\nrdiffParams \times \initialValues} \colon \R^\nrdiffParams \times \initialValues \to \EndValues
$
with $\nrdiffParams \in \N$ trainable parameters.
For any base model parameters 
$
	\approxParams \in \R^\nrbaseParams
$
with 
$
	\error(\adannBase_{\approxParams}) \in (0, \infty)
$
we then suggest training the difference model by minimizing the loss
$\differenceintLoss{\approxParams} \colon  \R^\nrdiffParams \to [0,\infty]$
given\footnote{We note that, in applications, the base model error 
$
	\error(\adannBase_{\approxParams}) 
= 
	\pr*{
		\Exp*{
			\lossmetric{\adannBase_{\approxParams}(\initialRV) - \solOp(\initialRV)}^2
		}
	}^{1/2}
$
has to be approximated by a suitable Monte Carlo method (cf.\ \cref{validation_error:eq1} in \cref{adanns_pseudocode}).} 
for all
$\theta \in \R^\nrdiffParams$
by
\begin{equation}
	\label{difference_training:eq3}
	\begin{split} 
		\differenceintLoss{\approxParams}(\theta)
		=
		\Exp*{
			\lossmetric{
				\diffann_\theta(\initialRV)
				-
					(\error(\adannBase_{\approxParams}))^{-1}
					\pr*{
						\solOp(\initialRV)
						-
						\adannBase_{\approxParams}(\initialRV)
					}
			}^2	
		}
	\end{split}
\end{equation}
by means of \SGD-type processes initialized in a standard random manner.
Note that \eqref{difference_training:eq3} suggests that, roughly speaking,
we train the difference model to approximate a scaled difference between the operator $\solOp$ and the base model $\adannBase_{\approxParams}$ instead of the difference itself.
This is a heuristic done so that the training samples for the difference model are roughly of the same scale irrespective of the approximation quality of the base model.
In our numerical simulations we found that this heuristic is of particular importance to ensure a successful training of the difference model when the base model achieves a very high approximation quality and thus the difference between the operator $\solOp$ and the base model $\adannBase_{\approxParams}$ is very small.

Combining the base model of \cref{sect:adann_base} and the difference model of this subsection, we define the \emph{full \ADANN\ model}
$
	\adann = (\adann_{(W, \theta, \epsilon)}(\IVvariable))_{(W, \theta, \epsilon, \IVvariable) \in \R^\nrbaseParams \times \R^\nrdiffParams \times \R \times \initialValues}
	\colon
	\R^\nrbaseParams \times \R^\nrdiffParams \times \R \times \initialValues \to \EndValues
$
by imposing
for all
$W \in \R^\nrbaseParams$,
$\theta \in \R^\nrdiffParams$,
$\epsilon \in \R$,
$\IVvariable \in \initialValues$
that
\begin{equation}
	\label{difference_training:eq4}
	\begin{split} 
		\adann_{(W, \theta, \epsilon)} (\IVvariable)
		=
		\adannBase_W(\IVvariable)
		+
		\epsilon
		\diffann_\theta (\IVvariable).
	\end{split}
\end{equation}
Observe that the scaling of the base model error in \eqref{difference_training:eq3} implies that 
for every 
$\approxParams \in \R^\nrbaseParams$,
$\theta \in \R^\nrdiffParams$
with
$
	\error(\adannBase_{\approxParams}) \in (0, \infty)
$
the 
$L^2$-error of the full \ADANN\ model $\adann_{(W, \theta, \error(\adannBase_{\approxParams}))}$ can be written as the product of the base and difference loss in the sense that
\begin{equation}
\label{difference_training:eq5}
\begin{split}
	\br*{
		\error(\adann_{(W, \theta, \error(\adannBase_{\approxParams}))})
	}^2
&=
		\Exp[\big]{
			\lossmetric{
				\adannBase_W(\IVvariable)
				+
				\error(\adannBase_{\approxParams})
				\diffann_\theta (\IVvariable)
				- 
				\solOp(\initialRV)
			}^2
		}
	\\
&=
	\br*{
		\error(\adannBase_{\approxParams})
	}^2
		\Exp[\big]{
			\lossmetric{
				\diffann_\theta (\IVvariable)
				+
				\pr*{
					\error(\adannBase_{\approxParams})
				}^{-1}
				\pr{
					\adannBase_W(\IVvariable)
					-
					\solOp(\initialRV)
				}
			}^2
		}
	\\
&=
	\br{\intLoss(\approxParams)}
	\br{\differenceintLoss{\approxParams}(\theta)}.
\end{split}
\end{equation}

\subsection{Optimization over base model initializations}
\label{sect:random_init}

\newcommand{\paramRV}{\mathfrak{p}}
\newcommand{\nrbaseSGDsteps}{N_{\text{Base}}}
\newcommand{\nrdiffSGDsteps}{N_{\text{Diff}}}
\newcommand{\optrun}{r^\ast}
\newcommand{\bboptimizer}{\mathfrak{o}}
\newcommand{\bbobj}{F}

The results of \SGD-type methods can strongly depend on their initializations and \eqref{base_model:eq4} suggests many choices of good initializations for the training of the base model.
In view of this, we propose to repeat the training of the full \ADANN\ model $\nrruns \in \N$ times
where in each training run the base model is initialized with different parameters out of the set $\{\initParams_p \colon p \in \params\}$ and the difference model is initialized in a standard random manner.
The best trained full \ADANN\ model over all training runs is then selected as the final approximation of the operator $\solOp$.

We aim to choose a sequence 
of initializations of the base model
which minimizes the error of the best full \ADANN\ model across all training runs.
Very roughly speaking, this can be considered to be a black box optimization (also called derivative free optimization) problem with noise (cf., e.g., \cite{Conn2009})
over the set of parameters $\params$
with the objective function 
given for all
$p \in \params$
by the (random) error of a trained full \ADANN\ model whose base model has been initialized with
$\initParams_p$
and whose difference model has been initialized in a standard random manner.
We abstractly model a black box optimization method for this optimization problem by a map
\begin{equation}
\label{random_init:eq1}
\begin{split}
	\bboptimizer \colon \pr*{ \cup_{r=0}^{\nrruns-1} (\params \times \R)^r } \times \Omega \to \params.
\end{split}
\end{equation}
Intuitively speaking, for every
	$r \in \{0, 1, \ldots, \nrruns-1\}$,
	$(p_1, e_1), \ldots, (p_r, e_r) \in \params \times [0, \infty)$
we think of 
$
	\initParams_{\bboptimizer((p_1, e_1), \ldots, (p_r, e_r))}
$ 
as the initialization 
for the training of the base model
in the (k+1)-th training run
 chosen by the black box optimization method 
$\bboptimizer$
 given that
$
	(p_1, e_1), \ldots, (p_k, e_k)
$
are the initialization parameters and the corresponding errors of the trained full \ADANN\ model of the first $k$ training runs.

For more concrete specifications of the black box optimization method $\bboptimizer$ we refer to \cref{sect:bbo} where we assume $\bboptimizer$ to be a grid-based black box optimization method (cf.\ \cref{sect:grid_optimizer}) and a heuristic exploration-exploitation black box optimization method (cf.\ \cref{sect:custom_optimizer}).
These two black box optimization approaches are the ones used in our numerical simulations in \cref{sect:simul}.

\subsection{Pseudocode description of the ADANN methodology}
\label{sect:adann_pseudocode}

In this section we describe in \cref{adanns_pseudocode} the \ADANN\ methodology in the form of a  pseudocode. 
Moreover, we also introduce a version of the \ADANN\ methodology in \cref{adanns_pseudocode_nodiff} which only relies on the base model of \cref{sect:adann_base} and the optimization approach of \cref{sect:random_init}
and does not include the difference model of \cref{sect:difference_training}.

In addition to the mathematical setting developed so far in \cref{sect:rough_overview},
we consider the element-wise base loss
$
	\loss 
	\colon
	\R^\nrbaseParams \times \initialValues 
	\to 
	[0,\infty)
$
given for all  
$W \in \R^\nrbaseParams$, 
$\IVvariable \in \initialValues$
by
\begin{equation}
\label{base_model_loss}
\begin{split}
	\loss(W, \IVvariable)
	=
	\lossmetric{\adannBase_{W}(\IVvariable) - \solOp(\IVvariable)}^2
\end{split}
\end{equation}
and the element-wise difference loss
$
	\diffLoss 
	\colon
	\R^\nrbaseParams \times \R^\nrdiffParams \times (0,\infty) \times \initialValues 
	\to 
	[0,\infty)
$
given for all
$W \in \R^\nrbaseParams$,
$\theta \in \R^\nrdiffParams$,
$\epsilon \in (0,\infty)$,
$\IVvariable \in \initialValues$
by
\begin{equation}
\label{diff_model_loss}
\begin{split}
	\diffLoss(W, \theta, \epsilon, \IVvariable)
	=
	\lossmetric*{
		\diffann_\theta(\IVvariable)
		-
		\pr*{
			\tfrac{1}{\epsilon}\pr*{
				\solOp(\IVvariable)
				-
				\adannBase_W(\IVvariable)
			}
		}
	}^2.
\end{split}
\end{equation}
Furthermore, we assume that for all 
$W \in \R^\nrbaseParams$,
$\epsilon \in (0,\infty)$,
$\IVvariable \in \initialValues$
the functions
$
	\loss(\cdot, \IVvariable) \colon \R^\nrbaseParams \to [0,\infty)
$
and
$
	\diffLoss(W, \cdot, \epsilon, \IVvariable) \colon \R^\nrdiffParams \to [0,\infty)
$
are differentiable.
With these new objects we are now in a position to formulate the \ADANN\ methodology in \cref{adanns_pseudocode} below.

\newcommand{\inlineComment}[1]{\textcolor{gray}{\textit{(#1)}}}
\newcommand{\mycomment}[1]{\hfill \textcolor{gray}{\textit{\# #1}}}
\newcommand{\multilinecomment}[3]{
\hfill 
\begin{minipage}[t]{#1}{
	\vspace{-#2}
	\begin{flushright}
		\textcolor{gray}{\textit{\# #3}}
	\end{flushright}
}
\end{minipage}}

\newcommand{\nrValSamples}{V}

\begin{myalgorithm}{\ADANN\ methodology}{adanns_pseudocode}
	\begin{algorithmic}[1] %
		\Statex\textbf{Setting:} All mathematical objects introduced in \cref{sect:rough_overview} above
		\Statex\textbf{Input:} 
		$\nrValSamples \in \N$ \inlineComment{Number of validation samples},
		$\gamma_{\text{Base}}, \gamma_{\text{Diff}} \in (0,\infty)$ \inlineComment{learning rates},
		$\batchsize_{\text{Base}}, \batchsize_{\text{Diff}} \in \N$ \inlineComment{batch sizes},
		$\nrbaseSGDsteps, \nrdiffSGDsteps \in \N$ \inlineComment{number of train steps}

		\Statex\textbf{Output:} 
		Approximation of $\solOp$

		\vspace{-2mm} %
		\noindent\hspace*{-7.7mm}\rule{\dimexpr\linewidth+8.8mm}{0.4pt}
		\vspace{-3mm} %

		\State 
		$\initialRV^{(\text{validate})}_{1}, \ldots, \initialRV^{(\text{validate})}_{\nrValSamples}
		\leftarrow \text{generate i.i.d.\ realizations of } \initialRV$
		\mycomment{Get validation samples}
		\For{$r = 1, \ldots, \nrruns$}
		
			\State 
			$
				\paramRV_r 
			\leftarrow
				\bboptimizer\pr[\big]{
					(\paramRV_1, \EndLoss_1), \ldots, (\paramRV_{r-1}, \EndLoss_{r-1})
				}
			$
			\mycomment{Choose initialization parameters}
			\State 
			$\baseSGD_r \leftarrow \initParams_{\paramRV_r}$
			\mycomment{Initialize base model}
			\For{$n = 1, \ldots, \nrbaseSGDsteps$}
				\State
				$\initialRV^{(\text{Base})}_{1}, \ldots, \initialRV^{(\text{Base})}_{\batchsize_{\text{Base}}} \leftarrow \text{generate i.i.d.\ realizations of } \initialRV$
				\label{training_samples:eq1}
				\mycomment{Get train samples}
				\State 
				$
					\baseSGD_r 
				\leftarrow
					\baseSGD_r
					- 
					\frac{\gamma_{\text{Base}}}{\batchsize_{\text{Base}}}
					\br[\big]{
						\sum_{b = 1}^{\batchsize_{\text{Base}}}
						(\nabla_W \loss)(\baseSGD_r, \initialRV^{(\text{Base})}_{b})
					}
				$
				\mycomment{Base model train step}
			\EndFor

			\State
			$
				\estLoss_r 
			\leftarrow
				\pr[\big]{
					\frac{1}{\nrValSamples} \br[\big]{
						\sum_{v = 1}^{\nrValSamples}
						\lossmetric{\adannBase_{\baseSGD_r}(\initialRV^{(\text{validate})}_{v} ) - \solOp(\initialRV^{(\text{validate})}_v)}^2
					}
				}^{\!1/2}	
			$
			\label{validation_error:eq1}
			\mycomment{Base model validation error}
			
			\State$\diffSGD_r \leftarrow \text{standard random initialization}$
			\mycomment{Initialize difference model}

			\For{$n = 1, \ldots, \nrdiffSGDsteps$}
				\State
				$\initialRV^{(\text{Diff})}_{1}, \ldots, \initialRV^{(\text{Diff})}_{\batchsize_{\text{Diff}}}
				\leftarrow \text{generate i.i.d.\ realizations of } \initialRV$
				\label{training_samples:eq2}
				\mycomment{Get train samples}

				\State 
				$
					\diffSGD_r
				\leftarrow
					\diffSGD_r
					- 
					\frac{\gamma_{\text{Diff}}}{\batchsize_{\text{Diff}}}
					\br[\big]{
						\sum_{b = 1}^{\batchsize_{\text{Diff}}}
						(\nabla_\theta \diffLoss)(
						\baseSGD_r, 
						\diffSGD_r, 
						\estLoss_r,
						\initialRV^{(\text{Diff})}_{b}
						)
					}
				$
				\mycomment{Difference model train step}

			\EndFor

			\State
			$
				\EndLoss_r
			\leftarrow
				\pr[\big]{
					\frac{1}{\nrValSamples} \br[\big]{
						\sum_{v = 1}^{\nrValSamples}
						\lossmetric{\adann_{\pr*{
								\baseSGD_r,
								\diffSGD_r,
								\estLoss_r
							}}(\initialRV^{(\text{validate})}_{v} ) - \solOp(\initialRV^{(\text{validate})}_v)}^2
					}
				}^{\!1/2}
			$
			\label{validation_error:eq2}
			\multilinecomment{4cm}{0.5cm}{Full \ADANN\ model  validation error}

		\EndFor

		\State 
		$\optrun \leftarrow \argmin_{r \in \{1, 2, \ldots, \nrruns\}} \EndLoss_r$
		\mycomment{Select full \ADANN\ model with lowest validation error}
		\State 
		\Return $\adann_{\pr*{
				\baseSGD^{(\optrun)},
				\diffSGD^{(\optrun)},
				\estLoss^{(\optrun)}
			}}$
	\end{algorithmic}
\end{myalgorithm}

Intuitively speaking, in some situations, adding a difference model to a base model might only introduce additional complexity without improving the approximation quality.
For this reason we also introduce a version of the \ADANN\ methodology whithout difference model in \cref{adanns_pseudocode_nodiff} below.
This version of the \ADANN\ methodology only relies on the base model of \cref{sect:adann_base} and the optimization approach of \cref{sect:random_init}.

\begin{myalgorithm}{\ADANN\ methodology without difference model}{adanns_pseudocode_nodiff}
	\begin{algorithmic}[1] %
		\Statex\textbf{Setting:} All mathematical objects introduced in \cref{sect:rough_overview} above
		\Statex\textbf{Input:}
		$\nrValSamples \in \N$ \inlineComment{Number of validation samples},
		$\gamma_{\text{Base}} \in (0,\infty)$ \inlineComment{learning rate},
		$\batchsize_{\text{Base}} \in \N$ \inlineComment{batch size},
		$\nrbaseSGDsteps \in \N$ \inlineComment{number of train steps}

		\Statex\textbf{Output:} 
		Approximation of $\solOp$
		
		\vspace{-2mm} %
		\noindent\hspace*{-7.7mm}\rule{\dimexpr\linewidth+8.8mm}{0.4pt}
		\vspace{-3mm} %
		
		\State 
		$\initialRV^{(\text{validate})}_{1}, \ldots, \initialRV^{(\text{validate})}_{\nrValSamples}
		\leftarrow \text{generate i.i.d.\ realizations of } \initialRV$
		\mycomment{Get validation samples}
		\For{$r = 1, \ldots, \nrruns$}
		
			\State 
			$
				\paramRV_r 
			\leftarrow
				\bboptimizer\pr[\big]{
					(\paramRV_1, \EndLoss_1), \ldots, (\paramRV_{r-1}, \EndLoss_{r-1})
				}
			$
			\mycomment{Choose initialization parameters}
			\State 
			$\baseSGD_r \leftarrow \initParams_{\paramRV_r}$
			\mycomment{Initialize base model}
			\For{$n = 1, \ldots, \nrbaseSGDsteps$}
				\State 
				$\initialRV^{(\text{Base})}_{1}, \ldots, \initialRV^{(\text{Base})}_{\batchsize_{\text{Base}}} \leftarrow \text{generate i.i.d.\ realizations of } \initialRV$
				\label{training_samples:eq3}
				\mycomment{Get train samples}
				\State 
				$
					\baseSGD_r 
				\leftarrow
					\baseSGD_r
					- 
					\frac{\gamma_{\text{Base}}}{\batchsize_{\text{Base}}}
					\br[\big]{
						\sum_{b = 1}^{\batchsize_{\text{Base}}}
						(\nabla_W \loss)(\baseSGD_r, \initialRV^{(\text{Base})}_{b})
					}
				$
				\mycomment{Base model train step}
			\EndFor

			\State
			$
				\EndLoss_r
			\leftarrow
				\pr[\big]{
					\frac{1}{\nrValSamples} \br[\big]{
						\sum_{v = 1}^{\nrValSamples}
						\lossmetric{
							\adannBase_{\baseSGD_r}(\initialRV^{(\text{validate})}_{v} ) - \solOp(\initialRV^{(\text{validate})}_v)}^2
					}
				}^{\!1/2}
			$
			\label{validation_error:eq3}
			\mycomment{Base model validation error}

		\EndFor

		\State 
		$\optrun \leftarrow \argmin_{r \in \{1, 2, \ldots, \nrruns\}} \EndLoss_r$
		\mycomment{Select base model with lowest validation error}
		\State 
		\Return 
		$
			\adannBase_{\baseSGD^{(\optrun)}}
		$
	\end{algorithmic}
\end{myalgorithm}

In the procedures described in \cref{adanns_pseudocode,adanns_pseudocode_nodiff}
above we made two major simplifications when compared to the methodology used in our numerical simulations in \cref{sect:simul}.
First, in our numerical simulations we train the models with the Adam optimizer with adaptive learning rates
(see \cref{sect:training} for more details on the training process in our numerical simulations).
However, for simplicity in \cref{adanns_pseudocode,adanns_pseudocode_nodiff} we only described the case of the plain-vanilla \SGD\ method with constant learning rates.
Second, in all our numerical simulations the operator $\solOp$ is not exactly known.
Therefore, in our numerical simulations, we replace the exact operator $\solOp$ with approximations of it 
	in the base loss (cf.\ \eqref{base_model_loss} above), 
	in the difference loss (cf.\ \eqref{diff_model_loss} above), 
	and 
	in the validation errors (cf.\ \cref{validation_error:eq1,validation_error:eq2} in \cref{adanns_pseudocode} and \cref{validation_error:eq3} in \cref{adanns_pseudocode_nodiff} above).
These approximations are computed ahead of the training process for all train and validation samples.
We note that, consequently, in our numerical simulations, the generation of training samples in 
\cref{training_samples:eq1,training_samples:eq2} in \cref{adanns_pseudocode} and
\cref{training_samples:eq3} in \cref{adanns_pseudocode_nodiff}
is done by drawing samples from a fixed training set.

\section{Derivation of a base model for semilinear heat PDEs}
\label{sect:semilinear_heat}

In this section we describe one way to design and initialize base models 
for the problem of approximating an operator mapping initial values to terminal values of a semilinear heat \PDE\ (cf.\ \cref{sect:semilinear_heat_setting}).
Specifically,
	we derive a family of approximation algorithms for the considered operator (cf.\ \cref{sect:semilinear_heat_design})
and 
	construct 
	an \ANN-type model together with a family of initializations for the model 
such that,
at these initializations, the model emulates the approximation algorithms (cf.\ \cref{sect:semilinear_heat_base}).
We will use very similar base models and initializations to apply the \ADANN\ methodology to a one-dimensional Sine-Gordon-type equation in our numerical simulations in \cref{sect:SineGordon1d}.

\subsection{One-dimensional semilinear heat PDEs}
\label{sect:semilinear_heat_setting}

\newcommand{\IVfunction}{g}
\newcommand{\sol}{u}
\newcommand{\nonlin}{f}
\newcommand{\Domain}{D}

\providecommand{\IVSpace}{}
\providecommand{\EndSpace}{}

\begingroup 
\providecommandordefault{\Domain}{(0,1)}
\providecommandordefault{\IVSpace}{C^2_{\text{per}}(\Domain, \R)}
\providecommandordefault{\EndSpace}{C^2_{\text{per}}(\Domain, \R)}

We now introduce the setting for the approximation problem considered in this section.
For this, assume the mathematical setting developed in \cref{sect:rough_overview},
let
$T \in (0,\infty)$,
let $\nonlin \colon \R \to \R$ be globally Lipschitz continuous,
and for every 
$\IVfunction \in \IVSpace$
let 
$\sol_\IVfunction \colon [0,T] \to \EndSpace$
be a mild solution
(cf., e.g., \cite{Cazenave1998,Pazy1983,Sell2002,Engel2000})
of the \PDE\
\begin{equation}
	\label{semilinear_heat:eq1}
	\begin{split} 
		\pr*{ \tfrac{\partial}{\partial t} \sol_\IVfunction } (t, x)
		=
		(\Delta_x \sol_\IVfunction)(t, x) + \nonlin (\sol_\IVfunction (t, x)),
	\qquad
		(t, x) \in [0,T] \times \Domain,
	\qquad
		\sol_\IVfunction(0) = \IVfunction
	\end{split}
\end{equation}
with periodic boundary conditions,
assume
$\initialValues = \EndValues = \IVSpace$,
and assume that the operator $\solOp \colon \initialValues \to \EndValues$ we want to approximate is
given for all
$\IVfunction \in \initialValues$ by
\begin{equation}
\label{semilinear_heat:eq1.1}
\begin{split}
	\solOp(\IVfunction)
	=
	\sol_\IVfunction(T).
\end{split}
\end{equation}

\subsection{Designing algorithms for the base model}
\label{sect:semilinear_heat_design}

We now derive a family of approximation algorithms 
for the operator $\solOp$ in \cref{semilinear_heat:eq1.1} which will serve as designing algorithms for the base model derived in this section.
Roughly speaking, the algorithms are based on discretizing the space domain of the \PDE\ in \cref{semilinear_heat:eq1} with the finite difference method and on discretizing the time domain of the \PDE\ in \cref{semilinear_heat:eq1} with a family of second order \LIRK\ methods.

\subsubsection{Spatial finite difference discretization}
\newcommand{\operator}{A}
\newcommand{\nrspacediscr}{N}
\newcommand{\odeSol}{\mathbf{u}}
\newcommand{\evalGrid}{\mathbf{e}}

For the spatial discretization of the \PDE\ in \eqref{semilinear_heat:eq1} consider
$\nrspacediscr \in \N$ grid points 
$\fx_1, \fx_2, \ldots, \fx_{\nrspacediscr} \in [0, 1]$ given for all
$i \in \{1, 2, \ldots,\nrspacediscr\}$ by $\fx_i = \frac{i-1}{\nrspacediscr}$,
let 
$\evalGrid \colon \initialValues \to \R^{\nrspacediscr}$ be the evaluation on the grid points given for all
$\IVfunction \in \initialValues$ by
\begin{equation}
	\label{semilinear_heat:eq2}
	\begin{split} 
		\evalGrid(\IVfunction)
		=
		(
			\IVfunction(\fx_1),
			\IVfunction(\fx_2),
			\ldots,
			\IVfunction(\fx_{\nrspacediscr})
		),
	\end{split}
\end{equation}
and consider the corresponding finite difference discretization of the Laplace operator on $(0,1)$
with periodic boundary conditions given by
\begin{equation}
	\label{semilinear_heat:eq3}
	\begin{split} 
		\operator
		=
		\nrspacediscr^2
		\begin{pmatrix}
			-2 & 1  & 0  & 0 & \cdots & 0 & 0  & 1    \\
			1  & -2 & 1  & 0 & \cdots & 0 & 0  & 0    \\
			0  & 1  & -2 & 1 & \cdots & 0 & 0  & 0    \\
			   &    &    &   & \ddots &   &    &      \\
			0  & 0  & 0  & 0 & \cdots & 1 & -2 & 1    \\ 
			1  & 0  & 0  & 0 & \cdots & 0 & 1  & -2 
		\end{pmatrix}
		\in
		\R^{\nrspacediscr \times \nrspacediscr}.
	\end{split}
\end{equation}
Under suitable assumptions, we expect for all 
$\IVfunction \in \initialValues$
that 
\begin{equation}
	\label{T_B_D}
	\begin{split} 
		\operator(
			\evalGrid(\IVfunction)
		)
	\approx
		\evalGrid(
				\Delta \IVfunction
		).
	\end{split}
\end{equation}
Using this spatial discretization on the \PDE\ in \eqref{semilinear_heat:eq1} results in an initial value \ODE.
Formally, for every
$\fg \in \R^\nrspacediscr$
let
$\odeSol_\fg \in C^{1}([0,T], \R^{\nrspacediscr})$ satisfy\footnote{Throughout this paper for every
	$h \colon \R \to \R$,
	$n \in \N$,
	$x = (x_1, x_2, \ldots x_n) \in \R^n$
	we denote by $h(x) \in \R^n$ the vector given by 
	$
		h(x) = (h(x_1), h(x_2), \ldots h(x_n))
	$.}
for all
$t \in [0,T]$
that
\begin{equation}
	\label{semilinear_heat:eq4}
	\begin{split} 
		\pr*{ \tfrac{\partial}{\partial t} \odeSol_\fg } (t)
		=
		\operator \odeSol_\fg (t) + \multdim{\nonlin}{\nrspacediscr} (\odeSol_\fg (t))
		\qandq
		\odeSol_\fg(0) = \fg.
	\end{split}
\end{equation}
Under suitable assumptions we expect for all
$\IVfunction \in \initialValues$ that
\begin{equation}
	\label{semilinear_heat:eq4.1}
	\begin{split} 
		\odeSol_{\evalGrid(\IVfunction)}(T)
	\approx
		\evalGrid(\sol_\IVfunction(T, \cdot))
	= 
		\evalGrid(\solOp(\IVfunction))
	.
	\end{split}
\end{equation}

\subsubsection{Temporal linearly implicit Runge-Kutta (LIRK) discretizations}
\label{sect:param_LIRK}
\newcommand{\maxstep}{\varepsilon}
\newcommand{\timestep}{H}
\newcommand{\LirkMethod}{\psi}
\newcommand{\interpolate}{\mathbf{i}}
\newcommand{\idmatrix}[1]{I}

In the next step we use a parametric family of second order \LIRK\ methods 
to discretize the \ODE\ in \eqref{semilinear_heat:eq4}.
We only introduce the family of \LIRK\ methods here and refer to \cref{sect:LIRK_derivation} for a more detailed derivation.

Specifically, 
let  
$\nrtimesteps \in \N$
be the numbers of time steps,
let 
$\timestep  = T/\nrtimesteps$
be the corresponding time step size,
let $\idmatrix{\nrspacediscr} \in \R^{\nrspacediscr \times \nrspacediscr}$ be the identity matrix,
assume $\params \subseteq (0,\infty)^2$,
and for all
parameters $p = (p_1, p_2) \in \params$ 
let the \LIRK\ time step
$
	\phi_{p}
	\colon 
	\R^\nrspacediscr \to \R^\nrspacediscr
$
satisfy for all
$U, k_1, k_2 \in \R^\nrspacediscr$
with
\begin{equation}
	\label{semilinear_heat:eq5}
	\begin{split} 
		\textstyle
		k_1
		=
		(\idmatrix{\nrspacediscr} - \timestep p_2 \operator)^{-1}
		\pr*{
			\operator
			U
			+
			\multdim{\nonlin}{\nrspacediscr}(U)
		}
		\qand
	\end{split}
\end{equation}
\begin{equation}
	\label{semilinear_heat:eq6}
	\begin{split} 
		k_2
		=
		(\idmatrix{\nrspacediscr} - \timestep p_2 \operator)^{-1}
		\pr*{
			\operator
			(U + \timestep 2p_1(\tfrac{1}{2} - p_2) k_1)
			+
			\multdim{\nonlin}{\nrspacediscr}(U + \timestep p_1 k_1)
		}
	\end{split}
\end{equation}
that
\begin{equation}
	\label{semilinear_heat:eq7}
	\begin{split} 
		\textstyle
		&\phi_{p}(U)
		=
		U
		+
		\timestep
		\br*{
			(1-\tfrac{1}{2 p_1}) k_1
			+
			(\tfrac{1}{2 p_1})
			k_2
		}.
	\end{split}
\end{equation}
For every choice of parameters
$p \in \params$
the corresponding \LIRK\ algorithm 
$
	\LirkMethod_{p}
	\colon 
	\R^\nrspacediscr \to \R^\nrspacediscr
$
for the \ODE\ in \eqref{semilinear_heat:eq4}
is then given by
\begin{equation}
	\label{semilinear_heat:eq8}
	\begin{split} 
		\LirkMethod_{p}
		=
		\underbrace{\phi_{p} \circ \ldots \circ \phi_{p}}_{\nrtimesteps \text{-times}}.
	\end{split}
\end{equation}
Under suitable assumptions 
(cf.\ \cref{semilinear_heat:eq4.1} above)
we then expect 
that for all 
$p \in \params$,
$\IVfunction \in \initialValues$
we have that
\begin{equation}
	\label{semilinear_heat:eq8.1}
	\begin{split} 
		\LirkMethod_{p}\pr*{
			\evalGrid(\IVfunction)
		}
	\approx
		\odeSol_{
			\evalGrid(\IVfunction)
		}(T)
	\approx 
		\evalGrid(\solOp(\IVfunction))
		.
	\end{split}
\end{equation}
Finally, let
$\interpolate \colon \R^{\nrspacediscr} \to \EndValues$ be an interpolation operator on the grid points in the sense that for all
$\fg = (\fg_1, \fg_2, \ldots, \fg_{\nrspacediscr}) \in \R^{\nrspacediscr}$,
$i \in \{1, 2, \ldots, \nrspacediscr\}$
we have that
\begin{equation}
\label{semilinear_heat:eq8.2}
\begin{split}
	\br{\interpolate(\fg)}(\fx_i)
	=
	\fg_i
\end{split}
\end{equation}
and assume that for every
$p \in \params$
the designing algorithm
$\propAlg_p \colon \initialValues \to \EndValues$
is given by
\begin{equation}
\label{semilinear_heat:eq8.3}
\begin{split}
	\propAlg_p
	=
	\interpolate
	\circ
	\LirkMethod_{p}
	\circ
	\evalGrid.
\end{split}
\end{equation}
Combining \cref{semilinear_heat:eq8.1,semilinear_heat:eq8.2,semilinear_heat:eq8.3}
we expect, under suitable assumptions, that for all
$p \in \params$,
$\IVfunction \in \initialValues$
we have that
\begin{equation}
	\label{semilinear_heat:eq8.4}
	\begin{split} 
		\propAlg_p(\IVfunction)
	=
		\interpolate(
			\LirkMethod_{p}(
				\evalGrid(\IVfunction)
			)
		)
	\approx
		\interpolate(
			\evalGrid(\solOp(\IVfunction))
		)
	\approx
		\solOp.
	\end{split}
\end{equation}
For a more rigorous error analysis for the designing algorithms in \cref{semilinear_heat:eq8.4} we refer, e.g., to \cite{tadmor2012review,Hairer2010,Thomee2006,Jovanovic2014} and the references therein.

\subsubsection{A compact reformulation of the designing algorithms}

\newcommand{\LIRKapprox}{\mathcal{U}}
\newcommand{\initMatrix}{\mathbf{w}}

To make the designing algorithms of \cref{sect:param_LIRK} amenable to be written as realizations of an \ANN-type base model we now present a more compact reformulation of the \LIRK\ time steps in \eqref{semilinear_heat:eq8}.
To this end
for every 
$p = (p_1, p_2) \in \params$
let
$\initMatrix_{p} = (\initMatrix_{p, i})_{i \in \{1, 2, \ldots, 5\}} \in (\R^{\nrspacediscr \times \nrspacediscr})^5$ satisfy
\begin{equation}
	\label{semilinear_heat:eq9}
	\begin{split} 
		\initMatrix_{p, 1}
		=
		(\idmatrix{\nrspacediscr} - \timestep p_2 \operator)^{-1}
		(\idmatrix{\nrspacediscr} + \timestep (1 - p_2) \operator) 
		+
		\timestep^2 (\tfrac{1}{2} - p_2) 
		\br*{
			(\idmatrix{\nrspacediscr} - \timestep p_2 \operator)^{-1}
			\operator
		}^2,
	\end{split}
\end{equation}
\begin{equation}
	\label{semilinear_heat:eq9.1}
	\begin{split} 
		\initMatrix_{p, 2}
		=
		\timestep
		(1-\tfrac{1}{2 p_1}) 
		(\idmatrix{\nrspacediscr} - \timestep p_2 \operator)^{-1}
		+
		\timestep^2
		(\tfrac{1}{2} - p_2) 
		(\idmatrix{\nrspacediscr} - \timestep p_2 \operator)^{-1}
		\operator
		(\idmatrix{\nrspacediscr} - \timestep p_2 \operator)^{-1},
	\end{split}
\end{equation}
\begin{equation}
	\label{semilinear_heat:eq9.2}
	\begin{split} 
		\initMatrix_{p, 3}
		=
		\timestep
		(\tfrac{1}{2 p_1})
		(\idmatrix{\nrspacediscr} - \timestep p_2 \operator)^{-1},
		\qquad
		\initMatrix_{p, 4}
		=
		(\idmatrix{\nrspacediscr} - \timestep p_2 \operator)^{-1}
		(\idmatrix{\nrspacediscr} + \timestep (p_1 - p_2) \operator),
	\end{split}
\end{equation}
\begin{equation}
	\label{semilinear_heat:eq9.3}
	\begin{split} 
		\andq
		\initMatrix_{p, 5}
		=
		\timestep p_1 
		(\idmatrix{\nrspacediscr} - \timestep p_2 \operator)^{-1}.
	\end{split}
\end{equation}
Note that 
\cref{semilinear_heat:eq5,semilinear_heat:eq6,semilinear_heat:eq7,semilinear_heat:eq9,semilinear_heat:eq9.1,semilinear_heat:eq9.2,semilinear_heat:eq9.3}
imply that
for all
$p \in \params$,
$U \in \R^\nrspacediscr$
we have that
\begin{equation}
	\label{semilinear_heat:eq10}
	\begin{split} 
		\phi_{p}(U)
		=
		\initMatrix_{p, 1} U 
		+ 
		\initMatrix_{p, 2} \multdim{\nonlin}{\nrspacediscr} (U)
		+
		\initMatrix_{p, 3} \multdim{\nonlin}{\nrspacediscr} \pr[\big]{
			\initMatrix_{p, 4} U 
			+ 
			\initMatrix_{p, 5} \multdim{\nonlin}{\nrspacediscr} (U)
		}.
	\end{split}
\end{equation}

\subsection{Designing the base model and its initializations}
\label{sect:semilinear_heat_base}

Roughly speaking, we propose to design the base model by considering the matrices in \eqref{semilinear_heat:eq10} as trainable parameters 
(see \cref{fig:base_model_schematic} for an illustration of the base model)
and to use \eqref{semilinear_heat:eq9}--\eqref{semilinear_heat:eq9.3} to define the initialization parameters.
More precisely, 
assume that 
$
	\nrbaseParams = 5 \nrspacediscr^{2d} \nrtimesteps
$,
assume\footnote{
	In order to structure the parameters of the base model, we slightly abuse the notation and identify 
	$
		\R^{\nrbaseParams}
	\simeq
		((\R^{\nrspacediscr^d \times \nrspacediscr^d})^5)^\nrtimesteps
	$.
} 
that the base model 
$
	\adannBase
\colon 
	((\R^{\nrspacediscr \times \nrspacediscr})^5)^\nrtimesteps 
	\times 
	\initialValues
	\to 
	\EndValues
$
is given for all
$W = ((W_{m, i})_{i \in \{1, 2, \ldots, 5\}})_{m \in \{1, 2, \ldots, \nrtimesteps\}} \allowbreak \in ((\R^{\nrspacediscr \times \nrspacediscr})^5)^\nrtimesteps$,
$\IVfunction \in \initialValues$,
$U_0, U_1, \ldots, U_\nrtimesteps \in \R^\nrspacediscr$
with
$U_0 = \evalGrid(\IVfunction)$
and
$
\forall \, m \in \{1, 2, \ldots, \nrtimesteps\} 
\colon
$
\begin{equation}
	\label{semilinear_heat:eq11}
	U_m 
	=  
	W_{m, 1} U_{m-1}
	+ 
	W_{m, 2} \multdim{\nonlin}{\nrspacediscr} (U_{m-1})
	+
	W_{m, 3} \multdim{\nonlin}{\nrspacediscr} \pr[\big]{
		W_{m, 4} U_{m-1}
		+ 
		W_{m, 5} \multdim{\nonlin}{\nrspacediscr} (U_{m-1})
	}
\end{equation}
by
\begin{equation}
	\label{semilinear_heat:eq12}
	\begin{split} 
		\adannBase_{W}(\IVfunction)
		=
		\interpolate(U_\nrtimesteps),
	\end{split}
\end{equation}
and for every 
	$p \in \params$ 
assume that the initialization parameters
$\initParams_p \in ((\R^{\nrspacediscr \times \nrspacediscr})^5)^\nrtimesteps$ are given by 
\begin{equation}
\label{T_B_D}
\begin{split}
	\initParams_p = \pr{
		\underbrace{
			\initMatrix_{p}, \ldots, \initMatrix_{p}
		}_{\nrtimesteps\text{-times}}
	}.
\end{split}
\end{equation}
Note that 
\enum{
	\eqref{semilinear_heat:eq8};
	\eqref{semilinear_heat:eq8.3};
	\eqref{semilinear_heat:eq10};
	\eqref{semilinear_heat:eq11};
	\eqref{semilinear_heat:eq12};
}[demonstrate]
that for all
$p \in \params$
it holds that
\begin{equation}
	\label{semilinear_heat:eq13}
	\begin{split} 
		\adannBase_{\initParams_p}
	=
		\pr[\big]{
			\interpolate
		\circ
			(\underbrace{\phi_{p} \circ \ldots \circ \phi_{p}}_{\nrtimesteps \text{-times}})
		\circ
			\evalGrid
		}	
	=
		\pr*{
			\interpolate
		\circ
			\LirkMethod_{p}
		\circ
			\evalGrid
		}
	=
		\propAlg_p.
	\end{split}
\end{equation}
Roughly speaking, we thus have constructed a base model which is able to reproduce all the designing algorithms we derived in \cref{sect:semilinear_heat_design}.
Combining this with 
\eqref{semilinear_heat:eq8.4}
implies that, under suitable assumptions, 
we have for all 
$p \in \params$
that
\begin{equation}
	\label{semilinear_heat:eq14}
	\begin{split} 
		&\adannBase_{\initParams_p}
	\approx	
		\solOp
	.
	\end{split}
\end{equation}

\begin{figure}
	\centering
	\includegraphics[width=0.9\linewidth]{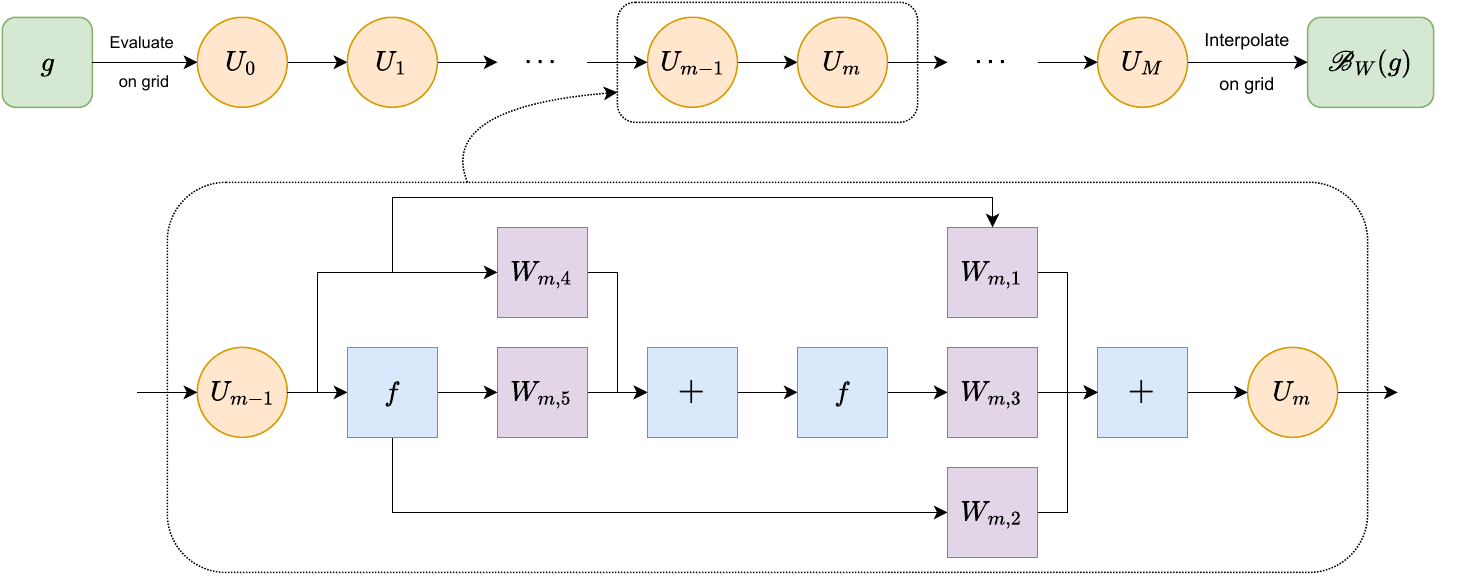}
	\caption{\label{fig:base_model_schematic}Graphical illustration for the base model defined in \eqref{semilinear_heat:eq11} and \eqref{semilinear_heat:eq12}.}
\end{figure}

\section{Numerical simulations}
\label{sect:simul}
\newcommand{\ann}{\mathscr{N}}
\newcommand{\laplaceFactor}{c}
\newcommand{\paramBoundary}{\mathfrak{p}}

In this section we numerically test the \ADANN\ methodology as described in \cref{sect:rough_overview}
in the case of four operators related to parametric \PDE\ problems
and compare its performance with existing operator learning architectures and classical numerical methods.
First, in \cref{sect:bbo} we describe the black box optimization methods used in the \ADANN\ methodology in our numerical simulations.
We then consider 
	the numerical approximation of
	operators mapping
	initial values to terminal values of one and two-dimensional Sine-Gordon-type equations in \cref{sect:SG_simul},
	the numerical approximation of
	an operator mapping initial values to terminal values of the viscous Burgers equation in \cref{sect:Burgers},
	and
	the numerical approximation of
	an operator mapping source terms to terminal values of a reaction-diffusion equation in \cref{sect:ReactionDiffusion}.

In every considered problem, all models (i.e., all base models, all difference models, and all other operator learning models)
are trained using the same training and validation set with the Adam optimizer with adaptive learning rates (cf.\ \cref{sect:training} for a detailed description of our adaptive training procedure).
Moreover, to provide a fair comparison between the \ADANN\ methodology and 
other operator learning architectures
we repeat the training of every operator learning model several times with different initializations and select the best performing trained model over all training runs as the approximation for that architecture.
In every considered problem, the $L^2$-errors of all methods are approximated using a Monte Carlo approximation based on the same test set.
The parameters chosen to generate train, validation, and test sets, as well as additional hyperparameters for each problem are listed in \cref{table:training_hyperparams} in \cref{sect:training}.

All the simulations were run on a remote machine on \texttt{https://vast.ai} equipped with an
\textsc{NVIDIA GeForce RTX 3090} GPU with 24 GB RAM
and an
\textsc{Xeon®} E5-2698 v3 CPU with 32 GB of total system RAM.
As the evaluation time of models on GPUs can be highly variable, we report the average evaluation time 
over $\nrEvalRuns$ test set evaluations for each method.
The code for all numerical simulations is available at 
\url{https://github.com/vwurstep/ADANNs}.

\subsection{Black box optimizers}
\label{sect:bbo}

One of the three main components of the \ADANN\ methodology described in \cref{sect:rough_overview} is to repeat the training of the full \ADANN\ model with different initializations and to use a black box optimization method over the highly specialized initializations of the base model to aim to minimize the error of the best trained full \ADANN\ model over all training runs (cf.\ \cref{sect:random_init}).
In principle any black box optimization method can be employed for the optimization problem arising in the \ADANN\ methodology and we tested several black box optimization methods from the scikit-opimize library (see \cite{Head2021}) such as Gaussian process optimization and decision tree based optimization, 
but they failed to give satisfactory results.
Consequently,
in our numerical simulations we restricted ourselves to two optimization approaches which we describe in detail in this section.

The first one is a simple grid-based black box optimization method. We describe this approach in \cref{sect:grid_optimizer}.
The second one is a heuristic black box optimization approach aiming to achieve an exploration-exploitation trade-off using an approximation of the objective function based on previously evaluated points.
We describe this approach in \cref{sect:custom_optimizer}.

\subsubsection{Grid-based black box optimizer}
\label{sect:grid_optimizer}

In this section we discuss a grid-based black box optimization approach to minimize a considered objective function.
Roughly speaking, a grid-based black box optimizer simply evaluates the objective function at all points on a grid irrespective of the outcome of previous evaluations.
To describe this approach more formally,
we will specify the black box optimizer 
$\bboptimizer$
introduced in the description of the \ADANN\ methodology in \cref{sect:random_init}
to such a grid-based black box optimizer.

For this, assume the setting developed in \cref{sect:rough_overview}, let 
$g_1, g_2, \ldots, g_\nrruns \in \params$,
and assume that for all
	$r \in \{0, 1, \ldots, \nrruns-1\}$,
	$(p_1, e_1), \ldots, (p_r, e_r) \in \params \times \R$
we have that
\begin{equation}
\label{T_B_D}
\begin{split}
	\bboptimizer((p_1, e_1), \ldots, (p_r, e_r))
=
	g_{r+1}.
\end{split}
\end{equation}
Roughly speaking, in this context 
we think of $\nrruns$ as the number of grid points and
we think of $g_1, g_2, \ldots, g_\nrruns$ as points on a grid in $\params$.

\subsubsection{Heuristic exploration-exploitation black box optimizer}
\label{sect:custom_optimizer}

\newcommand{\paramsDim}{\mathfrak{p}}
\newcommand{\nrRandomStarts}{Q}
\newcommand{\randomStartPoints}{\mathfrak{Q}}
\newcommand{\pastInfo}{\mathbf{p}}
\newcommand{\penaltyFactor}{c_1}
\newcommand{\inversionDecayRate}{c_2}
\newcommand{\interpolator}{\mathscr{r}}
\newcommand{\penalizedInterpolator}{\mathfrak{r}}
\newcommand{\unnormDensizy}{\mathfrak{D}}
\newcommand{\density}{\mathfrak{d}}
\newcommand{\bboRVs}{\mathfrak{X}}

In this section we discuss a heuristic exploration-exploitation black box optimization approach to minimize a considered objective function.
Roughly speaking, the proposed black box optimization approach can be divided into two phases.
In the first phase, the objective function is evaluated at several initial points, which may be chosen randomly, in its domain.
In the second phase, every new evaluation point is chosen by randomly sampling a point according to a probability density function which favors areas with lower approximate objective function values but penalizes areas which are close to previously evaluated points.
The approximate objective function values in the second phase are obtained
by a suitable regression or interpolation technique based on earlier evaluations of the objective function.
To describe this approach more formally,
we will specify the black box optimizer
$\bboptimizer$
introduced in the description of the \ADANN\ methodology in \cref{sect:random_init}
to such a heuristic exploration-exploitation black box optimizer.

For this
assume the setting developed in \cref{sect:rough_overview}, 
let 
	$\paramsDim \in \N$,
	$\penaltyFactor, \inversionDecayRate \in (0, \infty)$,
	$\nrRandomStarts \in \{0, 1, \ldots, \nrruns\}$,
assume that $\params$ is a compact subset of $\R^{\paramsDim}$,
let $\randomStartPoints_r \colon \Omega \to \params$, $r \in \{1, 2, \ldots, \nrRandomStarts\}$, be random variables,
let 
$
	\interpolator 
\colon 
	\pr*{ \cup_{r=0}^{\nrruns-1} (\params \times \R)^r } 
\to
	C(\params, (0,\infty))
$,
$
	\penalizedInterpolator
\colon
	\pr*{ \cup_{r=0}^{\nrruns-1} (\params \times \R)^r }
\to
	C(\params, (0,\infty))
$,
$
	\unnormDensizy
\colon
	\pr*{ \cup_{r=0}^{\nrruns-1} (\params \times \R)^r }
\to
	C(\params, (0,\infty))
$, 
and
$
	\density
\colon
	\pr*{ \cup_{r=0}^{\nrruns-1} (\params \times \R)^r }
\to
	C(\params, (0,\infty))
$
satisfy for all 
$r \in \{0, 1, \ldots, \nrruns-1\}$,
$\pastInfo = \allowbreak ((p_1, e_1), \allowbreak \ldots, (p_r, e_r)) \in (\params \times \R)^r$
that
\begin{equation}
\label{T_B_D}
\begin{split}
	\penalizedInterpolator_{\pastInfo}
=
	\interpolator_{\pastInfo}
	+
	\sum_{k=1}^{r}
		\tfrac{
			\penaltyFactor
			\min(\interpolator_{\pastInfo})
		}{
			\penaltyFactor + \norm{p_k - \,\cdot\,}_{\R^{\paramsDim}}
		},
\quad
	\unnormDensizy_{\pastInfo}
=
	\exp \pr*{
		- 
		\tfrac{
			\inversionDecayRate
			(\penalizedInterpolator_{\pastInfo} - \min(\interpolator_{\pastInfo}))
		}{
			\min(\interpolator_{\pastInfo})
		}
	},
\quad\text{and}\quad
	\density_{\pastInfo}
=
	\tfrac{
		\unnormDensizy_{\pastInfo}
	}{
		\int_{\params} \unnormDensizy_{\pastInfo}(q) d q
	},
\end{split}
\end{equation}
for every 
$\pastInfo \in \cup_{r=\nrRandomStarts}^{\nrruns-1} (\params \times \R)^r$
let 
$\bboRVs_{\pastInfo} \colon \Omega \to \params$ 
be a random variable which satisfies for all
$B \in \Borel(\params)$
that
\begin{equation}
\label{T_B_D}
\begin{split}
	\P \pr*{
		\bboRVs_{\pastInfo} \in B
	}
=
	\int_{B} \density_{\pastInfo}(q) d q,
\end{split}
\end{equation}
and assume for all
$r \in \{0, 1, \ldots, \nrruns-1\}$,
$\pastInfo = ((p_1, e_1), \ldots, (p_r, e_r)) \in (\params \times \R)^r$
that
\begin{equation}
\label{T_B_D}
\begin{split}
	\bboptimizer(\pastInfo)
=
	\begin{cases}
		\randomStartPoints_{r+1}, & r < \nrRandomStarts\\
		\bboRVs_{\pastInfo}, & r \geq \nrRandomStarts.
	\end{cases}
\end{split}
\end{equation}
Roughly speaking, 
	we think of $\randomStartPoints_1, \ldots, \randomStartPoints_\nrRandomStarts$ as the initial points at which the objective function is evaluated
and
for every number 
	$r \in \{\nrRandomStarts + 1, \nrRandomStarts + 2, \ldots, \nrruns\}$
and every sequence 
	$\pastInfo \in (\params \times \R)^r$ 
of previous evaluations 
\begin{enumerate}[label=(\roman*)]
	\item 
	we think of $\bboRVs_{\pastInfo}$ as the point at which the objective function is evaluated in the $r$-th step of the optimization process given that the previous evaluations are given by $\pastInfo$,

	\item 
	we think of $\density_{\pastInfo}$ as the probability density function obtained from normalizing $\unnormDensizy_{\pastInfo}$
	which is used to sample the point $\bboRVs_{\pastInfo}$,

	\item	
	we think of $\unnormDensizy_{\pastInfo}$ as an unnormalized probability density function which inverses the values of $\penalizedInterpolator_{\pastInfo}$, that is, $\unnormDensizy_{\pastInfo}$ assigns higher probabilities to areas where $\penalizedInterpolator_{\pastInfo}$ has lower values and lower probabilities where $\penalizedInterpolator_{\pastInfo}$ has higher values,

	\item 
	we think of $\penalizedInterpolator_{\pastInfo}$ as a penalized version of $\interpolator_{\pastInfo}$ where the function values close to previous evaluation points in $\pastInfo$ are increased,
	and

	\item
	we think of $\interpolator_{\pastInfo}$ as an approximation of the objective function based on the previous evaluations in $\pastInfo$.
\end{enumerate}
In all our numerical simulations we choose the hyperparameters 
$\penaltyFactor = 0.005$ and $\inversionDecayRate = 100$ 
and 
we take $\interpolator$ to be a smoothened RBF interpolation as implemented in \cite{SciPy2021}.

\subsection{Sine-Gordon-type equation}
\label{sect:SG_simul}
\newcommand{\EndVariable}{h}

\begingroup
\providecommandordefault{\Domain}{(0,\SGonedSpaceSize)^d}
\providecommandordefault{\IVSpace}{H^2_{\text{per}}(\Domain; \R)}
\providecommandordefault{\EndSpace}{H^2_{\text{per}}(\Domain; \R)}

In this section we test the \ADANN\ methodology without difference model as described in \cref{sect:adann_pseudocode} (cf.\ \cref{adanns_pseudocode_nodiff})
in the case of operators mapping initial values to terminal values of one and two-dimensional Sine-Gordon-type equations.
We introduce below the considered Sine-Gordon-type equations and the corresponding operators.
The base models used for both the one and two-dimensional case are defined in \cref{sect:SG_base} and the results of our numerical simulations are presented in \cref{sect:SineGordon1d} for the one-dimensional case and in \cref{sect:SineGordon2d} for the two-dimensional case.

Throughout \cref{sect:SG_simul} assume the mathematical setting developed in \cref{sect:rough_overview},
let
$d \in \{1, 2\}$,
$T = \SGonedT$,
$\laplaceFactor = \SGonedLaplaceFactor$,
for every 
$\IVfunction \in \IVSpace$
let 
$\sol_\IVfunction \colon [0,T] \to \EndSpace$
be a mild solution of the \PDE\
\begin{equation}
	\label{SG:eq1}
	\begin{split} 
		\pr*{ \tfrac{\partial}{\partial t} \sol_\IVfunction } (t, x)
		=
		\laplaceFactor
		(\Delta_x \sol_\IVfunction)(t, x) + \sin (\sol_\IVfunction (t, x)),
	\qquad 
		(t, x) \in [0,T] \times \Domain,
	\qquad
		\sol_\IVfunction(0) = \IVfunction
	\end{split}
\end{equation} 
with periodic boundary conditions,
assume 
$\initialValues = \EndValues = \IVSpace$,
and 
assume 
that the operator
$\solOp \colon \initialValues \to \EndValues$
we want to approximate is
given for all
$\IVfunction \in \initialValues$ by
\begin{equation}
\label{SG:eq2}
\begin{split}
	\solOp(\IVfunction)
	=
	\sol_\IVfunction(T).
\end{split}
\end{equation}
Moreover, assume that the initial value 
$\initialRV \colon \Omega \to \initialValues$
is
$\SGonedInitDistr$-distributed
where $\Delta_x$ is the Laplace operator on $\Ltwo{\Domain}$
with periodic boundary conditions,
fix a space discretization $\nrspacediscr \in \N$,
and assume for all
$\EndVariable \in \EndValues$
that
\begin{equation}
\label{SG:eq3}
\begin{split}
		\lossmetric{\EndVariable}^2
	=
		\frac{1}{\nrspacediscr^d}
		\br*{
			\sum_{\fx \in \{\frac{0}{\nrspacediscr}, \frac{1}{\nrspacediscr}, \ldots, \frac{\nrspacediscr-1}{\nrspacediscr}\}^d}
				(\EndVariable(\fx))^2
		}
	\approx
		\int_{\Domain} 
			(\EndVariable(x))^2
		d x
		.
\end{split}
\end{equation}
We recall that our goal is to find an approximation $\solOpAlt \colon \initialValues \to \EndValues$ of the operator in \eqref{SG:eq2} which minimizes the $L^2$-error
$
	\error(\solOpAlt)
\in 
	[0, \infty] 
$ given by
\begin{equation}
	\label{SG:eq4}
	\begin{split}
		\error(\solOpAlt)
	=
		\pr[\big]{
			\Exp[\big]{
				\lossmetric{\solOpAlt(\initialRV) - \solOp(\initialRV)}^2
			}
		}^{1/2}
	\approx
		\pr*{
			\Exp*{
				\int_{\Domain} 
					\pr[\big]{
						\solOpAlt(\initialRV)(x) - 	\sol_\initialRV(T, x)
					}^2
				d x
			}
		}^{1/2}.
	\end{split}
\end{equation}

\subsubsection{Base model for the Sine-Gordon-type equation}
\label{sect:SG_base}

We now describe the base model that we use in the \ADANN\ methodology to approximate the operator in \eqref{SG:eq2}.
Roughly speaking, we use the base model derived in \cref{sect:semilinear_heat} with the initialization parameters adjusted to the considered Sine-Gordon-type equations.
A graphical illustration for this base model can be found in \cref{fig:base_model_schematic}.

Specifically, 
let
$\nrtimesteps \in \N$,
let $\evalGrid \colon \initialValues \to \R^{\nrspacediscr^d}$ be an evaluation operator on the  grid $\{\frac{0}{\nrspacediscr}, \frac{1}{\nrspacediscr}, \ldots, \frac{\nrspacediscr-1}{\nrspacediscr}\}^d$ 
and
let $\interpolate \colon \R^{\nrspacediscr^d} \to \EndValues$ be a corresponding interpolation operator (cf.\ \cref{semilinear_heat:eq2,semilinear_heat:eq8.2}),
let $\nonlin = \sin$, and
assume that 
$\nrbaseParams = 5 \nrspacediscr^{2d} \nrtimesteps$.
We then assume\footnote{
	In order to structure the parameters of the base model, we slightly abuse the notation and identify 
	$
		\R^{\nrbaseParams}
	\simeq
		((\R^{\nrspacediscr^d \times \nrspacediscr^d})^5)^\nrtimesteps
	$.
}  
that the base model
$
	\adannBase
	\colon 
	((\R^{\nrspacediscr^d \times \nrspacediscr^d})^5)^\nrtimesteps 
	\times 
	\initialValues
	\to 
	\EndValues
$
is given for all 
$W = ((W_{m, i})_{i \in \{1, 2, \ldots, 5\}})_{m \in \{1, 2, \ldots, \nrtimesteps\}} \allowbreak \in ((\R^{\nrspacediscr^d \times \nrspacediscr^d})^5)^\nrtimesteps$,
$\IVfunction \in \initialValues$,
$U_0, U_1, \ldots, U_\nrtimesteps \in \R^{\nrspacediscr^d}$
with
$U_0 = \evalGrid(\IVfunction) $
and
$
\forall \, m \in \{1, 2, \ldots, \nrtimesteps\} 
\colon
$
\begin{equation}
	\label{SG:eq5}
	U_m 
	=  
	W_{m, 1} U_{m-1}
	+ 
	W_{m, 2} \multdim{\nonlin}{\nrspacediscr^d} (U_{m-1})
	+
	W_{m, 3} \multdim{\nonlin}{\nrspacediscr^d} \pr[\big]{
		W_{m, 4} U_{m-1}
		+ 
		W_{m, 5} \multdim{\nonlin}{\nrspacediscr^d} (U_{m-1})
	}
\end{equation}
by
\begin{equation}
	\label{SG:eq6}
	\begin{split} 
		\adannBase_{W}(\IVfunction)
		=
		\interpolate(U_\nrtimesteps).
	\end{split}
\end{equation}

To define the family of initialization parameters 
$\initParams_p \in ((\R^{\nrspacediscr^d \times \nrspacediscr^d})^5)^\nrtimesteps$, $p \in \params$,
for the base model
let  
$\timestep = T/\nrtimesteps$,
let $\idmatrix{\nrspacediscr^d} \in \R^{\nrspacediscr^d \times \nrspacediscr^d}$ be the identity matrix,
let $\operator \in \R^{\nrspacediscr^d \times \nrspacediscr^d}$ be the finite difference discretization of the Laplace operator on $\Domain$ with periodic boundary conditions,
corresponding to the evaluation operator $\evalGrid$ 
(cf. \eqref{semilinear_heat:eq3} for a definition of $\operator$ in the case $d = 1$),
assume $\params = \SGonedParams$,
and
for every 
$p = (p_1, p_2) \in \params$
let
$\initMatrix_{p} = (\initMatrix_{p, i})_{i \in \{1, 2, \ldots, 5\}} \in (\R^{\nrspacediscr^d \times \nrspacediscr^d})^5$
satisfy
\begin{equation}
	\label{SG:eq7}
	\begin{split} 
		\initMatrix_{p, 1}
		=
		(\idmatrix{\nrspacediscr} - \timestep p_2 \laplaceFactor \operator)^{-1}
		(\idmatrix{\nrspacediscr} + \timestep (1 - p_2) \laplaceFactor \operator) 
		+
		\timestep^2 (\tfrac{1}{2} - p_2) 
		\br*{
			(\idmatrix{\nrspacediscr} - \timestep p_2 \laplaceFactor \operator)^{-1}
			\laplaceFactor \operator
		}^2,
	\end{split}
\end{equation}
\begin{equation}
	\label{SG:eq7.1}
	\begin{split} 
		\initMatrix_{p, 2}
		=
		\timestep
		(1-\tfrac{1}{2 p_1}) 
		(\idmatrix{\nrspacediscr} - \timestep p_2 \laplaceFactor \operator)^{-1}
		+
		\timestep^2
		(\tfrac{1}{2} - p_2) 
		(\idmatrix{\nrspacediscr} - \timestep p_2 \laplaceFactor \operator)^{-1}
		\laplaceFactor \operator
		(\idmatrix{\nrspacediscr} - \timestep p_2 \laplaceFactor \operator)^{-1},
	\end{split}
\end{equation}
\begin{equation}
	\label{SG:eq7.2}
	\begin{split} 
		\initMatrix_{p, 3}
		=
		\timestep
		(\tfrac{1}{2 p_1})
		(\idmatrix{\nrspacediscr} - \timestep p_2 \laplaceFactor \operator)^{-1},
		\qquad
		\initMatrix_{p, 4}
		=
		(\idmatrix{\nrspacediscr} - \timestep p_2 \laplaceFactor \operator)^{-1}
		(\idmatrix{\nrspacediscr} + \timestep (p_1 - p_2) \laplaceFactor \operator),
	\end{split}
\end{equation}
\begin{equation}
	\label{SG:eq7.3}
	\begin{split} 
	\andq
		\initMatrix_{p, 5}
	=
		\timestep p_1 
		(\idmatrix{\nrspacediscr} - \timestep p_2 \laplaceFactor \operator)^{-1}.
	\end{split}
\end{equation}
We then assume that
for every 
$p \in \params$
the parameters
$\initParams_p \in ((\R^{\nrspacediscr^d \times \nrspacediscr^d})^5)^\nrtimesteps$
are given by
\begin{equation}
\label{SG:eq7.4}
\begin{split}
	\initParams_p = \pr{
		\underbrace{
			\initMatrix_{p}, \ldots, \initMatrix_{p}
		}_{\nrtimesteps\text{-times}}
	}.
\end{split}
\end{equation}
Roughly speaking, for every $p \in \params$ we have that 
\begin{equation}
	\label{SG:eq8}
	\begin{split} 
		\adannBase_{\initParams_p}
	\approx
		\solOp
	\end{split}
\end{equation}
corresponds to an approximation of 
the Sine-Gordon-type equation in \eqref{SG:eq1}
based on a finite difference discretization in space and a \LIRK\ approximation in time where the parameters $p$ corresponds to the parameters of the \LIRK\ method (cf. \cref{sect:LIRK_derivation}).

\subsubsection{Numerical results for the one-dimensional Sine-Gordon-type equation}
\label{sect:SineGordon1d}

\begingroup
\providecommandordefault{\eqname}{Semilinear_heat_1-dimensional_T_2.0_space_size_1.0_laplace_factor_0.01_nonlin_Sine_var_100000_decay_rate_2_offset_316.22776601683796_inner_decay_1.0}

In this section we present numerical results for the approximation of the operator in \eqref{SG:eq1} in the case $d = 1$.
We test the \ADANN\ methodology without difference model (cf.\ \cref{adanns_pseudocode_nodiff})
with
	the base model and the corresponding initializations defined in \cref{sect:SG_base},
	parameter space $\params = \SGonedParams$, 
	space discretization $\nrspacediscr = \SGonedSpaceStep$,
	number of time steps $\nrtimesteps \in \{2, 4, 8\}$,
	and
	both
	a grid-based black box optimizer as described in \cref{sect:grid_optimizer}
	(see rows 10-12 in \cref{table:SineGordon1d} and \cref{fig:SineGordon1d_grid})
	and 
	our heuristic exploration-exploitation black box optimizer as described in \cref{sect:custom_optimizer}
	(see rows 13-15 in \cref{table:SineGordon1d} and \cref{fig:SineGordon1d_opt}).
We also test different \ANN\ models with \GELU\ activation function
(see rows 1-3 in \cref{table:SineGordon1d}),
\FNO\ models 
(see rows 4-6 in \cref{table:SineGordon1d}),
and classical methods 
(see rows 7-9 in \cref{table:SineGordon1d})
for comparison.
As classical methods we use the untrained base model
$
	\adannBase_{\initParams_{(0.5, 0.5)}}
$ 
with $\nrtimesteps \in \{2, 4, 8\}$ time steps, 
corresponding, roughly speaking, 
to a finite difference discretization in space and a Crank-Nicolson explicit midpoint \LIRK\ discretization in time (cf.\ \cref{sect:crank_nicolson}).
The performance of all considered methods is summarized in \cref{table:SineGordon1d} and graphically illustrated in \cref{fig:SineGordon1d}.
In addition, some approximations for a randomly chosen test sample are shown in \cref{fig:SineGordon1d_sample}.

\begin{table} 
	\tiny
	\resizebox{\textwidth}{!}{
		\csvreader[
			tabular=|c|c|c|c|c|,
			separator=semicolon,
			table head=
			\hline 
			\thead{Method} &  
			\thead{Estimated \\ $L^2$-error \\ in \cref{SG:eq4}} &
			\thead{Average evaluation time \\ for $\SGonedNrTestSamples$ test samples \\ over $\nrEvalRuns$ runs (in s)} &
			\thead{Number \\ of trainable \\ parameters} &
			\thead{Precomputation\\time (in s)}
			\\\hline,
			late after line=\\\hline
		]
		{1_numbers/rounded_methods_data_\eqname.csv}
		{
			Method=\method, 
			L2_error = \llerror, 
			nr_params = \numparams, 
			training_time = \traintime, 
			test_time = \evaltime
		}
		{\method& \llerror& \evaltime &\numparams&\traintime}%
	}
	\caption{\label{table:SineGordon1d} 
	Comparison of the performance of different methods for the approximation of the operator in \eqref{SG:eq2} mapping initial values to terminal values of the Sine-Gordon-type equation in \cref{SG:eq1} in the case $d = 1$.}
\end{table}

\begin{figure}
	\centering
	\includegraphics[width=0.7\linewidth]{0_plots/1dsemilin/error_scatter_plot_\eqname.pdf}
	\caption{\label{fig:SineGordon1d}Graphical illustration of the performance of the methods in \cref{table:SineGordon1d}.}
\end{figure}

\begin{figure}
	
	\includegraphics[width=\linewidth]{0_plots/1dsemilin/Results_grid/grid_error_overview_\eqname.pdf}
	
	\caption{\label{fig:SineGordon1d_grid}
		Illustration of the \ADANN\ methodology without difference model (cf.\ \cref{adanns_pseudocode_nodiff})
		with a grid-based black box optimizer applied to 
		the approximation of the operator in \eqref{SG:eq2} based on the Sine-Gordon-type equation in \cref{SG:eq1} in the case $d=1$.
		\emph{Left}: Test errors of the base models prior to training as a function of the parameters used for initialization.
		\emph{Right}: Test errors of the trained base models as a function of the parameters used for initialization.
	}
\end{figure}

\begin{figure}
	\includegraphics[width=\linewidth]{0_plots/1dsemilin/Results_opt/opt_scatter_combined_\eqname.pdf}
	
	\caption{\label{fig:SineGordon1d_opt} 
		Illustration of the 
		\ADANN\ methodology without difference model (cf.\ \cref{adanns_pseudocode_nodiff})
		with our heuristic exploration-exploitation black box optimizer 
		applied to the approximation of the operator in \eqref{SG:eq2} 
		mapping initial values to terminal values of the Sine-Gordon-type equation in \cref{SG:eq1} in the case $d=1$.
		\emph{Left}: Test errors of trained base models as a function of the parameters used for initialization.  Increasing scatter sizes indicate higher training run numbers.
		\emph{Middle}: The same test errors represented in the order in which they appeared in the black box optimization process.
		\emph{Right}: Coordinates of the chosen parameters in the black box optimization process.
	}
\end{figure}

\begin{figure} 
	\minipage{0.48\textwidth}
	\expandafter\includegraphics\expandafter[width=\linewidth]{"0_plots/1dsemilin/Sample_plots/ol_plots_\eqname_0.pdf"}
	\endminipage
	\minipage{0.48\textwidth}
	\expandafter\includegraphics\expandafter[width=\linewidth]{"0_plots/1dsemilin/Sample_plots/adann_classical_plots_\eqname_0.pdf"}
	\endminipage
	\caption{\label{fig:SineGordon1d_sample}
		Example approximation plots for a randomly chosen sample from the test set for the Sine-Gordon-type equation in \cref{SG:eq1} in the case $d=1$.
		\emph{Left}: \ANN\ and \FNO\ approximations.
		\emph{Right}: Classical and \ADANN\ approximations.
	}
\end{figure}

\endgroup

\clearpage
\subsubsection{Numerical results for the two-dimensional Sine-Gordon-type equation}
\label{sect:SineGordon2d}

\begingroup
\providecommandordefault{\eqname}{Semilinear_heat_2-dimensional_T_2.0_space_size_1.0_laplace_factor_0.01_nonlin_Sine_var_100000_decay_rate_2_offset_316.22776601683796_inner_decay_1.0}

In this section we present numerical results for 
the approximation of the operator in \eqref{SG:eq2} in the case $d = 2$.
We test the \ADANN\ methodology without difference model (cf.\ \cref{adanns_pseudocode_nodiff})
with
	the base model and the corresponding initializations defined in \cref{sect:Burgers_base},
	space discretization $\nrspacediscr = \SGtwodSpaceStep$,
	number of time steps $\nrtimesteps \in \{2, 4, 8\}$,
	and
	our heuristic exploration-exploitation black box optimizer as described in \cref{sect:custom_optimizer}
	(see rows 10-12 in \cref{table:SineGordon2d} and \cref{fig:SineGordon2d_opt}).
We also test different \ANN\ models with \GELU\ activation function
(see rows 1-3 in \cref{table:SineGordon2d}),
\FNO\ models 
(see rows 4-6 in \cref{table:SineGordon2d}),
and classical methods 
(see rows 7-9 in \cref{table:SineGordon2d})
for comparison.
As classical methods we use the untrained base model
$
	\adannBase_{\initParams_{(0.5, 0.5)}}
$ 
with $\nrtimesteps \in \{2, 4, 8\}$ time steps, 
corresponding, roughly speaking, 
to a finite difference discretization in space and a Crank-Nicolson explicit midpoint \LIRK\ discretization in time (cf.\ \cref{sect:crank_nicolson}).
The performance of all considered methods are summarized in \cref{table:SineGordon2d} and graphically illustrated in \cref{fig:SineGordon2d}.
In addition, some approximations for a randomly chosen test sample are shown in \cref{fig:SineGordon2d_sample}.

\begin{table} 
	\tiny
	\resizebox{\textwidth}{!}{
		\csvreader[
			tabular=|c|c|c|c|c|,
			separator=semicolon,
			table head=
			\hline 
			\thead{Method} &  
			\thead{Estimated \\ $L^2$-error \\ in \cref{SG:eq4}} &
			\thead{Average evaluation time \\ for $\SGtwodNrTestSamples$ test samples \\ over $\nrEvalRuns$ runs (in s)} &
			\thead{Number \\ of trainable \\ parameters} &
			\thead{Precomputation\\time (in s)}
			\\\hline,
			late after line=\\\hline
		]
		{1_numbers/rounded_methods_data_\eqname.csv}
		{
			Method=\method, 
			L2_error = \llerror, 
			nr_params = \numparams, 
			training_time = \traintime, 
			test_time = \evaltime
		}
		{\method& \llerror& \evaltime &\numparams&\traintime}%
	}
	\caption{\label{table:SineGordon2d}
	Comparison of the performance of different methods for the approximation of the operator in \eqref{SG:eq2} 
	mapping initial values to terminal values of the Sine-Gordon-type equation in \cref{SG:eq1} in the case $d = 2$.
	}
\end{table}

\begin{figure}
	\centering
	\includegraphics[width=0.7\linewidth]{0_plots/2dsemilin/error_scatter_plot_\eqname.pdf}
	\caption{\label{fig:SineGordon2d}
		Graphical illustration of the performance of the methods in \cref{table:SineGordon2d}.
	}
\end{figure}

\begin{figure}
\centering
 \includegraphics[width=0.9\linewidth]{0_plots/2dsemilin/Results_opt/opt_scatter_combined_\eqname.pdf}
 	\caption{\label{fig:SineGordon2d_opt}Illustration of the 
	\ADANN\ methodology without difference model (cf.\ \cref{adanns_pseudocode_nodiff})
	with our heuristic exploration-exploitation black box optimizer
	applied to the approximation of the operator in \eqref{SG:eq2}
	mapping initial values to terminal values of the Sine-Gordon-type equation in \cref{SG:eq1} in the case $d=2$.
	 \emph{Left}: Test errors of trained base models as a function of parameters used for initialization.  Increasing scatter sizes indicate higher training run numbers.
	 \emph{Middle}: The same test errors represented in the order in which they appeared in the black box optimization process.
	 \emph{Right}: Coordinates of the chosen parameters in the black box optimization process.
 	}
\end{figure}

\begin{figure}
	\centering
	\includegraphics[width=0.7\linewidth]{0_plots/2dsemilin/Sample_plots/all_plots_\eqname_0.pdf}
	\caption{Example approximation plots for a randomly chosen sample from the test set for the Sine-Gordon-type equation in \cref{SG:eq1} in the case $d=2$.}
	\label{fig:SineGordon2d_sample}
\end{figure}

\endgroup

\endgroup

\clearpage
\subsection{Viscous Burgers equation}
\label{sect:Burgers}
\newcommand{\firstOderDiff}{E}

\begingroup
\providecommandordefault{\eqname}{Burgers_T1.0_S6.283185307179586_nu0.1_var1000_decay3.0_offset9.999999999999998_innerdecay2.0}

\providecommandordefault{\Domain}{(0, \BurgersSpaceSize)}
\providecommandordefault{\IVSpace}{H^2_\text{per}(\Domain; \R)}
\providecommandordefault{\EndSpace}{H^2_\text{per}(\Domain; \R)}

In this section we test the \ADANN\ methodology as described in \cref{sect:adann_pseudocode} (cf.\ \cref{adanns_pseudocode,adanns_pseudocode_nodiff})
in the case of an operator mapping initial values to terminal values of the viscous Burgers equation.
We introduce below the viscous Burgers equation in conservative form and the corresponding operator.
The base model for the \ADANN\ methodology is defined in \cref{sect:Burgers_base} and the results of our numerical simulations are presented in \cref{sect:Burgers_simul}.

Throughout \cref{sect:Burgers}
assume the mathematical setting developed in \cref{sect:rough_overview},
let
$T = \BurgersT$,
$\laplaceFactor = \BurgersLaplaceFactor$,
for every 
$\IVfunction \in \IVSpace$
let 
$\sol_\IVfunction \colon [0, T] \to \EndSpace$ 
be a mild solution of the \PDE\
\begin{equation}
\label{Burgers:eq1}
\begin{split}
	\pr*{ \tfrac{\partial}{\partial t} \sol_\IVfunction } (t, x)
	=
	\laplaceFactor
	(\Delta_x \sol_\IVfunction)(t, x) 
	- 
	\tfrac{1}{2}
	\pr*{\tfrac{\partial}{\partial t} \sol_\IVfunction^2}(t, x),
\quad
	(t, x) \in [0, T] \times \Domain,
\quad
	\sol_\IVfunction(0) = \IVfunction
\end{split}
\end{equation}
with periodic boundary conditions,
assume 
$\initialValues = \EndValues = \IVSpace$,
and assume that the operator
$\solOp \colon \initialValues \to \EndValues$
we want to approximate is
given for all
$\IVfunction \in \initialValues$ by
\begin{equation}
\label{Burgers:eq2}
\begin{split}
	\solOp(\IVfunction)
	=
	\sol_\IVfunction(T).
\end{split}
\end{equation}
Moreover, assume that the initial value 
$\initialRV \colon \Omega \to \initialValues$
is
$\BurgersInitDistr$-distributed
where $\Delta_x$ is the Laplace operator on $\Ltwo{\Domain}$
with periodic boundary conditions, 
fix a space discretization $\nrspacediscr=\BurgersSpaceStep$,
and assume for all
$\EndVariable \in \EndValues$
that
\begin{equation}
\label{Burgers:eq3}
\begin{split}
		\lossmetric{\EndVariable}^2
	=
		\frac{
			\BurgersSpaceSize
		}{\nrspacediscr}
		\br*{
			\sum_{\fx \in \{\frac{0}{\nrspacediscr}, \frac{1}{\nrspacediscr}, \ldots, \frac{\nrspacediscr-1}{\nrspacediscr}\}}
				(\EndVariable(\BurgersSpaceSize \fx))^2
		}
	\approx
		\int_0^{\BurgersSpaceSize}
			(\EndVariable(x))^2
		d x
		.
\end{split}
\end{equation}
Recall that our goal is to find an approximation $\solOpAlt \colon \initialValues \to \EndValues$ of the operator in \eqref{Burgers:eq2} which minimizes the $L^2$-error
$
	\error(\solOpAlt)
\in 
	[0, \infty] 
$ given by
\begin{equation}
	\label{Burgers:eq4}
	\begin{split}
		\error(\solOpAlt)
	=
		\pr[\big]{
			\Exp[\big]{
				\lossmetric{\solOpAlt(\initialRV) - \solOp(\initialRV)}^2
			}
		}^{1/2}
	\approx
		\pr*{
			\Exp*{
				\int_0^{\BurgersSpaceSize}
					\pr[\big]{
						\solOpAlt(\initialRV)(x) - \sol_\initialRV(T, x)
					}^2
				d x
			}
		}^{1/2}.
	\end{split}
\end{equation}

\subsubsection{Base model for the viscous Burgers equation}
\label{sect:Burgers_base}

We now describe the base model that we use in the \ADANN\ methodology to approximate the operator in \eqref{Burgers:eq2}.
Roughly speaking, we use the same base model architecture as in \cref{sect:semilinear_heat_base} but we additionally integrate the first order derivative operator coming from the nonlinearity into the learnable parameters of the model.
A graphical illustration for this base model can be found in \cref{fig:base_model_schematic}.

More precisely, let
$\nrtimesteps \in \N$,
let $\evalGrid \colon \initialValues \to \R^{\nrspacediscr}$ be the evaluation operator on the grid
$\{\frac{0}{\nrspacediscr}, \frac{\BurgersSpaceSize}{\nrspacediscr}, \ldots, \frac{\BurgersSpaceSize(\nrspacediscr-1)}{\nrspacediscr}\}$
and
let $\interpolate \colon \R^{\nrspacediscr} \to \EndValues$ be a corresponding interpolation operator
(cf.\ \cref{semilinear_heat:eq2,semilinear_heat:eq8.2}),
let $\nonlin \colon \R \to \R$ satisfy for all
$x \in \R$
that
$
	\nonlin(x) = -\tfrac{1}{2} x^2
$,
and
assume that 
$\nrbaseParams = 5 \nrspacediscr^2 \nrtimesteps$.
We then assume\footnote{
	In order to structure the parameters of the base model, we slightly abuse the notation and identify 
	$
		\R^{\nrbaseParams}
	\simeq
		((\R^{\nrspacediscr \times \nrspacediscr})^5)^\nrtimesteps
	$.
}
that the base model
$
	\adannBase
	\colon 
	((\R^{\nrspacediscr \times \nrspacediscr})^5)^\nrtimesteps 
	\times 
	\initialValues
	\to 
	\EndValues
$
is given for all
$W = ((W_{m, i})_{i \in \{1, 2, \ldots, 5\}})_{m \in \{1, 2, \ldots, \nrtimesteps\}} \allowbreak \in ((\R^{\nrspacediscr \times \nrspacediscr})^5)^\nrtimesteps$,
$\IVfunction \in \initialValues$,
$U_0, U_1, \ldots, U_\nrtimesteps \in \R^{\nrspacediscr}$
with
$U_0 = \evalGrid(\IVfunction) $
and
$
\forall \, m \in \{1, 2, \ldots, \nrtimesteps\} 
\colon
$
\begin{equation}
	\label{Burgers:eq5}
	U_m 
	=  
	W_{m, 1} U_{m-1}
	+ 
	W_{m, 2} \multdim{\nonlin}{\nrspacediscr} (U_{m-1})
	+
	W_{m, 3} \multdim{\nonlin}{\nrspacediscr} \pr[\big]{
		W_{m, 4} U_{m-1}
		+ 
		W_{m, 5} \multdim{\nonlin}{\nrspacediscr} (U_{m-1})
	}
\end{equation}
by
\begin{equation}
	\label{Burgers:eq6}
	\begin{split} 
		\adannBase_{W}(\IVfunction)
		=
		\interpolate(U_\nrtimesteps).
	\end{split}
\end{equation}

To define the family of initialization parameters 
$\initParams_p \in ((\R^{\nrspacediscr \times \nrspacediscr})^5)^\nrtimesteps$, $p \in \params$,
for the base model
let  
$\timestep = T/\nrtimesteps$,
let $\idmatrix{\nrspacediscr} \in \R^{\nrspacediscr \times \nrspacediscr}$ be the identity matrix,
let $\operator \in \R^{\nrspacediscr \times \nrspacediscr}$ be the finite difference discretization of the Laplace operator on $\Domain$ with periodic boundary conditions corresponding to the evaluation operator $\evalGrid$ (cf.\ \eqref{semilinear_heat:eq3} for the definition of $\operator$ for the domain $[0,1]$),
let $\firstOderDiff \in \R^{\nrspacediscr \times \nrspacediscr}$ be the finite difference discretization of the first order derivative operator with periodic boundary conditions on $\Domain$ 
given by 
\begin{equation}
\label{T_B_D}
\begin{split}
	\firstOderDiff
	=
	\frac{\BurgersSpaceSize}{2}
	\begin{pmatrix}
		0 & 1 & 0 & 0 &\cdots & 0 & 0 & -1 \\
		-1 & 0 & 1 & 0 &\cdots & 0 & 0 & 0 \\
		0 & -1 & 0 & 1 &\cdots & 0 & 0 & 0 \\
		&&&&\ddots& \\
		0 & 0 & 0 & 0 &\cdots & -1 & 0 & 1 \\
		1 & 0 & 0 & 0 &\cdots & 0 & -1 & 0 
	\end{pmatrix}
	\in 
	\R^{\nrspacediscr \times \nrspacediscr},
\end{split}
\end{equation}
assume $\params = \BurgersParams$,
and
for every 
$p = (p_1, p_2) \in \params$
let
$\initMatrix_{p} = (\initMatrix_{p, i})_{i \in \{1, 2, \ldots, 5\}} \in (\R^{\nrspacediscr \times \nrspacediscr})^5$
satisfy
\begin{equation}
	\label{Burgers:eq7}
	\begin{split} 
		\initMatrix_{p, 1}
		=
		(\idmatrix{\nrspacediscr} - \timestep p_2 \laplaceFactor \operator)^{-1}
		(\idmatrix{\nrspacediscr} + \timestep (1 - p_2) \laplaceFactor \operator) 
		+
		\timestep^2 (\tfrac{1}{2} - p_2) 
		\br*{
			(\idmatrix{\nrspacediscr} - \timestep p_2 \laplaceFactor \operator)^{-1}
			\laplaceFactor \operator
		}^2,
	\end{split}
\end{equation}
\begin{equation}
	\label{Burgers:eq7.1}
	\begin{split} 
			\initMatrix_{p, 2}
		=
			\br*{
				\timestep
				(1-\tfrac{1}{2 p_1}) 
				(\idmatrix{\nrspacediscr} - \timestep p_2 \laplaceFactor \operator)^{-1}
				+
				\timestep^2
				(\tfrac{1}{2} - p_2) 
				(\idmatrix{\nrspacediscr} - \timestep p_2 \laplaceFactor \operator)^{-1}
				\laplaceFactor \operator
				(\idmatrix{\nrspacediscr} - \timestep p_2 \laplaceFactor \operator)^{-1}
			}
			\firstOderDiff,
	\end{split}
\end{equation}
\begin{equation}
	\label{Burgers:eq7.2}
	\begin{split} 
			\initMatrix_{p, 3}
		=
			\timestep
			(\tfrac{1}{2 p_1})
			(\idmatrix{\nrspacediscr} - \timestep p_2 \laplaceFactor \operator)^{-1}
			\firstOderDiff,
		\qquad
			\initMatrix_{p, 4}
		=
			(\idmatrix{\nrspacediscr} - \timestep p_2 \laplaceFactor \operator)^{-1}
			(\idmatrix{\nrspacediscr} + \timestep (p_1 - p_2) \laplaceFactor \operator),
	\end{split}
\end{equation}
\begin{equation}
	\label{Burgers:eq7.3}
	\begin{split} 
	\andq
		\initMatrix_{p, 5}
	=
		\timestep p_1 
		(\idmatrix{\nrspacediscr} - \timestep p_2 \laplaceFactor \operator)^{-1}
		\firstOderDiff.
	\end{split}
\end{equation}
We then assume that
for every 
$p \in \params$
the parameters 
$\initParams_p \in ((\R^{\nrspacediscr^d \times \nrspacediscr^d})^5)^\nrtimesteps$
are given by
\begin{equation}
\label{T_B_D}
\begin{split}
	\initParams_p = \pr{
		\underbrace{
			\initMatrix_{p}, \ldots, \initMatrix_{p}
		}_{\nrtimesteps\text{-times}}
	}.
\end{split}
\end{equation}
Roughly speaking, for every $p \in \params$ we have that 
\begin{equation}
	\label{Burgers:eq8}
	\begin{split} 
		\adannBase_{\initParams_p}
	\approx
		\solOp
	\end{split}
\end{equation}
corresponds to an approximation of the viscous Burgers equation in \eqref{Burgers:eq1}
based on a finite difference discretization in space and a \LIRK\ discretization in time where the parameter $p$ corresponds to the parameters of the \LIRK\ method (cf. \cref{sect:LIRK_derivation}).

\subsubsection{Numerical results for the viscous Burgers equation}
\label{sect:Burgers_simul}

In this section we present numerical results for the approximation of the operator in \eqref{Burgers:eq2}.
We test the \ADANN\ methodology with 
(see rows 13--15 in \cref{table:Burgers} and \cref{fig:Burgers_grid}) 
and without difference model
(see rows 10--12 in \cref{table:Burgers} and \cref{fig:Burgers_grid})
with
	the base model and the corresponding initializations defined in \cref{sect:Burgers_base},
	parameter space $\params = \BurgersParams$,
	space discretization $\nrspacediscr = \BurgersSpaceStep$,
	number of time steps $\nrtimesteps \in \{2, 4, 8\}$,
	grid-based black box optimizer as described in \cref{sect:grid_optimizer},
	and
	difference model given by an \ANN\ with architecture $\BurgersDMarch$.
We also test different \ANN\ models with \GELU\ activation function
(see rows 1-3 in \cref{table:Burgers}),
\FNO\ models 
(see rows 4-6 in \cref{table:Burgers}),
and classical methods 
(see rows 7-9 in \cref{table:Burgers})
for comparison.
As classical methods we use the untrained base model
$
	\adannBase_{\initParams_{(0.5, 0.5)}}
$ 
with $\nrtimesteps \in \{2, 4, 8\}$ time steps, 
corresponding, roughly speaking, 
to a finite difference discretization in space and a Crank-Nicolson explicit midpoint \LIRK\ discretization in time (cf.\ \cref{sect:crank_nicolson}).
The performance of all considered methods is summarized in \cref{table:Burgers} and graphically illustrated in \cref{fig:Burgers}.
In addition, some approximations for a randomly chosen test sample are shown in \cref{fig:Burgers_sample}.

\begin{table} 
	\tiny
	\resizebox{\textwidth}{!}{
		\csvreader[
			tabular=|c|c|c|c|c|,
			separator=semicolon,
			table head=
			\hline 
			\thead{Method} &  
			\thead{Estimated \\ $L^2$-error \\ in \cref{Burgers:eq4}} &
			\thead{Average evaluation time \\ for $\BurgersNrTestSamples$ test samples \\ over $\nrEvalRuns$ runs (in s)} &
			\thead{Number \\ of trainable \\ parameters} &
			\thead{Precomputation\\time (in s)}
			\\\hline,
			late after line=\\\hline
		]
		{1_numbers/rounded_methods_data_\eqname.csv}
		{
			Method=\method, 
			L2_error = \llerror, 
			nr_params = \numparams, 
			training_time = \traintime, 
			test_time = \evaltime
		}
		{\method& \llerror& \evaltime &\numparams&\traintime}%
	}
	\caption{\label{table:Burgers}
	Comparison of the performance of different methods for the approximation of the operator in \eqref{Burgers:eq2}
	mapping initial values to terminal values of the viscous Burgers equation in \cref{Burgers:eq1}.
	}
\end{table}

\begin{figure}
	\centering
	\includegraphics[width=0.7\linewidth]{0_plots/Burgers/error_scatter_plot_\eqname.pdf}
	\caption{\label{fig:Burgers}
		Graphical illustration of the performance of the methods in \cref{table:Burgers}.
	}
\end{figure}

\begin{figure}
	\includegraphics[width=\linewidth]{0_plots/Burgers/Results_grid/grid_error_overview_\eqname.pdf}
	\caption{\label{fig:Burgers_grid}
		Illustration of the \ADANN\ methodology with and without difference model (cf.\ \cref{adanns_pseudocode,adanns_pseudocode_nodiff}) applied to the approximation of the operator in \eqref{Burgers:eq2} 
		mapping initial values to terminal values of the viscous Burgers equation in \cref{Burgers:eq1}.
		\emph{Left}: Test errors of the base models prior to training as a function of the parameters used for initialization.
		\emph{Middle}: Test errors of the trained base models as a function of the parameters used for initialization.
		\emph{Right}: Test errors of the trained full \ADANN\ models as a function of the parameters used for initialization of the base model.
	}
\end{figure}

\begin{figure}
	\minipage{0.48\textwidth}
	\expandafter\includegraphics\expandafter[width=\linewidth]{"0_plots/Burgers/Sample_plots/ol_plots_\eqname_0.pdf"}
	\endminipage
	\minipage{0.48\textwidth}
	\expandafter\includegraphics\expandafter[width=\linewidth]{"0_plots/Burgers/Sample_plots/ADANN_full_-_grid_plots_\eqname_0.pdf"}
	\endminipage
	\caption{\label{fig:Burgers_sample}
			Example approximation plots for a randomly chosen sample from the test set for the viscous Burgers equation in \cref{Burgers:eq1}.
			\emph{Left}: \ANN\ and \FNO\ approximations.
			\emph{Right}: \ADANN\ approximations.
	}
\end{figure}

\endgroup

\clearpage
\subsection{Reaction-diffusion equation}
\label{sect:ReactionDiffusion}

\begingroup
\providecommandordefault{\eqname}{ReactionDiffusion_T1.0_S2.0_nu0.05_k2.0_nonlinAllenCahn_var10000_decay2_offset100.0_innerdecay1.0}

\providecommandordefault{\reactionRate}{k}

\providecommandordefault{\Domain}{(0,\ReactDiffSpaceSize)}
\providecommandordefault{\IVSpace}{H^2_\text{per}(\Domain; \R)}
\providecommandordefault{\EndSpace}{H^2_\text{per}(\Domain; \R)}

In this section we apply the \ADANN\ methodology 
as described in \cref{sect:adann_pseudocode} (cf.\ \cref{adanns_pseudocode,adanns_pseudocode_nodiff})
in the case of an operator mapping source terms to terminal values of a reaction-diffusion equation.
The considered operator is inspired by the reaction-diffusion equation in \cite[Section 4.3]{Lu2021}.
We introduce below the considered reaction-diffusion equation and the corresponding operator.
The base model for the \ADANN\ methodology is defined in \cref{sect:RD_base} and the results of our numerical simulations are presented in \cref{sect:RD_simul}.

Throughout \cref{sect:ReactionDiffusion}
assume the mathematical setting developed in \cref{sect:rough_overview},
let
$T = \ReactDiffT$,
$\laplaceFactor = \ReactDiffLaplaceFactor$,
$\reactionRate = \ReactDiffReactionRate$,
for every 
$\IVfunction \in \IVSpace$
let 
$\sol_\IVfunction \colon [0, T] \to \EndSpace$
be a mild solution of the \PDE\
\begin{equation}
	\label{ReactionDiffusion:eq1}
	\begin{split} 
		\pr*{ \tfrac{\partial}{\partial t} \sol_\IVfunction } (t, x)
		=
		\laplaceFactor
		(\Delta_x \sol_\IVfunction)(t, x) 
		+
		\reactionRate 
		\pr*{
			\sol_\IVfunction(t, x) - (\sol_\IVfunction(t, x))^3
		}
		+
		\IVfunction(x),
	\end{split}
\end{equation}
\begin{equation}
\label{ReactionDiffusion:eq1.1}
	\begin{split}
		(t, x) \in [0, T] \times \Domain,
		\qquad
		\sol_\IVfunction(0) = 0
	\end{split}
\end{equation} 
with periodic boundary conditions,
assume
$\initialValues = \EndValues = \IVSpace$,
and assume that the operator
$\solOp \colon \initialValues \to \EndValues$
we want to approximate is
given for all
$\IVfunction \in \initialValues$ by
\begin{equation}
\label{ReactionDiffusion:eq2}
\begin{split}
	\solOp(\IVfunction)
	=
	\sol_\IVfunction(T).
\end{split}
\end{equation}
Moreover, assume that the source term
$\initialRV \colon \Omega \to \initialValues$
is
$\ReactDiffInitDistr$-distributed
where $\Delta_x$ is the Laplace operator on 
$\Ltwo{\Domain}$
with periodic boundary conditions,
fix a space discretization $\nrspacediscr=128$,
and assume for all
$\EndVariable \in \EndValues$
that
\begin{equation}
\label{ReactionDiffusion:eq3}
\begin{split}
		\lossmetric{\EndVariable}^2
	=
		\frac{1}{\nrspacediscr}
		\br*{
			\sum_{\fx \in \{\frac{0}{\nrspacediscr}, \frac{1}{\nrspacediscr}, \ldots, \frac{\nrspacediscr-1}{\nrspacediscr}\}}
				(\EndVariable(\fx))^2
		}
	\approx
		\int_0^{\ReactDiffSpaceSize}
			(\EndVariable(x))^2
		d x
		.
\end{split}
\end{equation}
Recall that our goal is to find an approximation $\solOpAlt \colon \initialValues \to \EndValues$ of the operator in \eqref{base_model:eq1} which minimizes the $L^2$-error
$
	\error(\solOpAlt)
\in 
	[0, \infty] 
$ given by
\begin{equation}
	\label{ReactionDiffusion:eq4}
	\begin{split}
		\error(\solOpAlt)
	=
		\pr[\big]{
			\Exp[\big]{
				\lossmetric{\solOpAlt(\initialRV) - \solOp(\initialRV)}^2
			}
		}^{1/2}
	\approx
		\pr*{
			\Exp*{
				\int_0^{\ReactDiffSpaceSize}
					\pr[\big]{
						\solOpAlt(\initialRV)(x) - \sol_\initialRV(T, x)
					}^2
				d x
			}
		}^{1/2}.
	\end{split}
\end{equation}

\subsubsection{Base model for the reaction-diffusion equation}
\label{sect:RD_base}
\newcommand{\sourceTerm}{\fg}

We now describe the base model that we use in the \ADANN\ methodology to approximate the operator in \eqref{ReactionDiffusion:eq2}.
Very roughly speaking, we use a similar approach as in \cref{sect:semilinear_heat_base} to define the base model with the adjustment that the initial value is fixed and the source term is involved in every layer through additional learnable parameters.
A graphical illustration for this base model can be found in \cref{fig:RD_base_model}.

More precisely, let
$\nrtimesteps \in \N$,
let $\evalGrid \colon \initialValues \to \R^{\nrspacediscr}$ be the evaluation operator on the  grid 
$ \{\frac{0}{\nrspacediscr}, \frac{\ReactDiffSpaceSize}{\nrspacediscr}, \ldots, \frac{(\nrspacediscr-1)}{\nrspacediscr}\}$ and
let $\interpolate \colon \R^{\nrspacediscr} \to \EndValues$ be a corresponding interpolation operator 
(cf.\ \cref{semilinear_heat:eq2,semilinear_heat:eq8.2}),
let $\nonlin \colon \R \to \R$ satisfy for all
$x \in \R$
that
$
	\nonlin(x) = \reactionRate (x - x^3)
$,
and
assume that 
$\nrbaseParams = 7\nrspacediscr^{2} \nrtimesteps$.
We then assume\footnote{
	In order to structure the parameters of the base model, we slightly abuse the notation and identify 
	$
		\R^{\nrbaseParams}
	\simeq
		((\R^{\nrspacediscr \times \nrspacediscr})^7)^\nrtimesteps
	$.
}
that the base model
$
	\adannBase
	\colon 
	((\R^{\nrspacediscr \times \nrspacediscr})^7)^\nrtimesteps 
	\times 
	\initialValues
	\to 
	\EndValues
$
satisfies for all
$W = ((W_{m, i})_{i \in \{1, 2, \ldots, 7\}})_{m \in \{1, 2, \ldots, \nrtimesteps\}} \allowbreak \in ((\R^{\nrspacediscr \times \nrspacediscr})^7)^\nrtimesteps$,
$\IVfunction \in \initialValues$,
$U_0, U_1, \ldots, U_\nrtimesteps \in \R^{\nrspacediscr}$
with
$U_0 = 0$,
$\sourceTerm = \evalGrid(\IVfunction)$,
and
$
\forall \, m \in \{1, 2, \ldots, \nrtimesteps\} 
\colon
$
\begin{equation}
	\label{ReactionDiffusion:eq5}
	U_m 
	=  
	W_{m, 1} U_{m-1}
	+ 
	W_{m, 2} \multdim{\nonlin}{\nrspacediscr} (U_{m-1})
	+
	W_{m, 6} \sourceTerm
	+
	W_{m, 3} \multdim{\nonlin}{\nrspacediscr} \pr[\big]{
		W_{m, 4} U_{m-1}
		+ 
		W_{m, 5} \multdim{\nonlin}{\nrspacediscr} (U_{m-1})
		+
		W_{m, 7} \sourceTerm
	}
\end{equation}
that
\begin{equation}
	\label{ReactionDiffusion:eq6}
	\begin{split} 
		\adannBase_{W}(\IVfunction)
		=
		\interpolate(U_\nrtimesteps).
	\end{split}
\end{equation}

\begin{figure}
	\centering
	\includegraphics[width=0.9\linewidth]{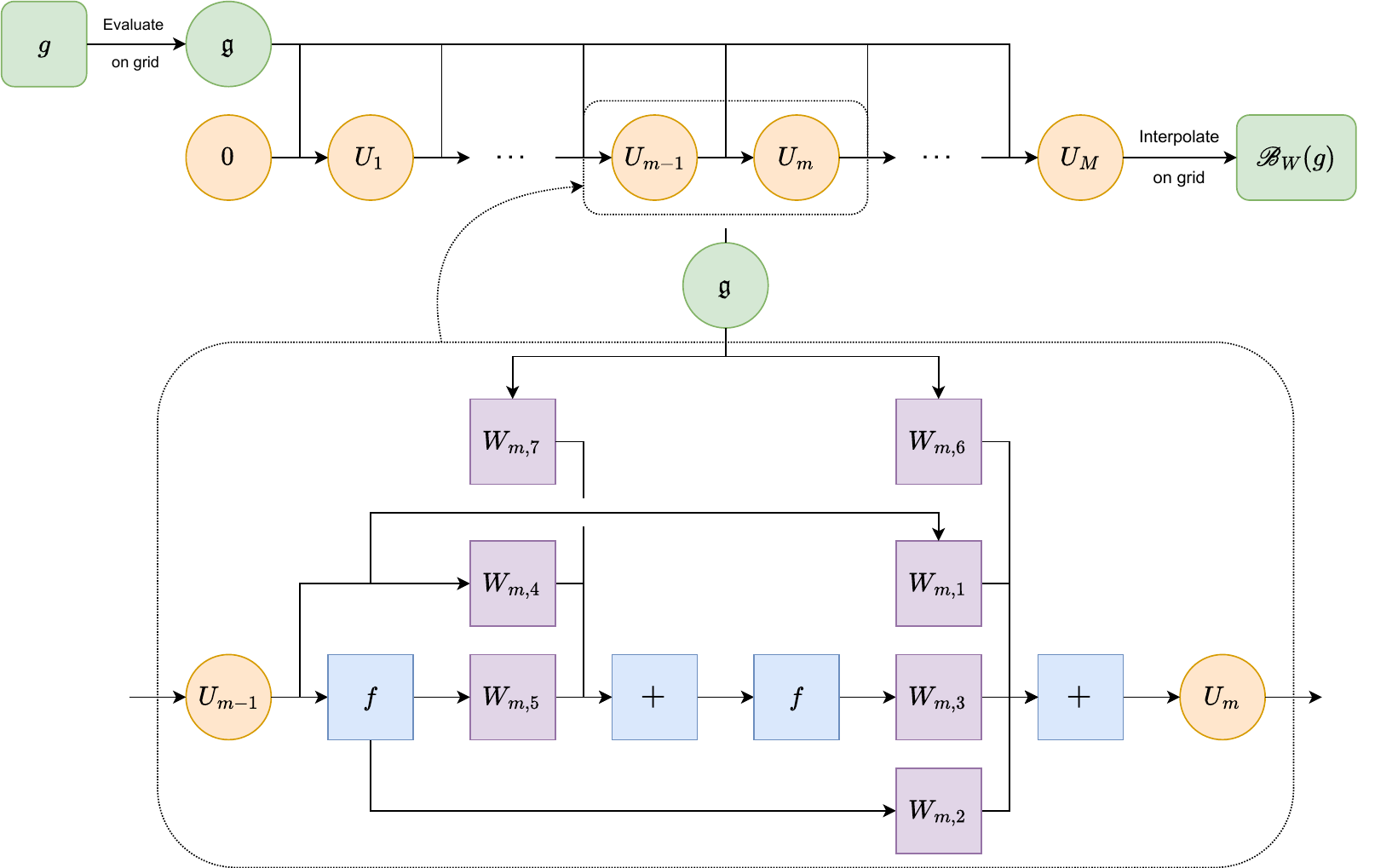}
	\caption{\label{fig:RD_base_model} Graphical illustration for the base model defined in \eqref{ReactionDiffusion:eq5} and \eqref{ReactionDiffusion:eq6}.}
\end{figure}

To define the family of initialization parameters
$\initParams_p \in ((\R^{\nrspacediscr \times \nrspacediscr})^7)^\nrtimesteps$, $p \in \params$,
for the base model
let  
$\timestep = T/\nrtimesteps$,
let $\idmatrix{\nrspacediscr} \in \R^{\nrspacediscr \times \nrspacediscr}$ be the identity matrix,
let $\operator \in \R^{\nrspacediscr \times \nrspacediscr}$ be the finite difference discretization of the Laplace operator on $\Domain$ with periodic boundary conditions 
corresponding to the evaluation operator $\evalGrid$ 
(cf.\ \eqref{semilinear_heat:eq3} for the definition of $\operator$),
assume $\params = \ReactDiffParams$,
and
for every 
$p = (p_1, p_2) \in \params$
let
$\initMatrix_{p} = (\initMatrix_{p, i})_{i \in \{1, 2, \ldots, 7\}} \in (\R^{\nrspacediscr \times \nrspacediscr})^7$
satisfy
\begin{equation}
	\label{ReactionDiffusion:eq7}
	\begin{split} 
		\initMatrix_{p, 1}
		=
		(\idmatrix{\nrspacediscr} - \timestep p_2 \laplaceFactor \operator)^{-1}
		(\idmatrix{\nrspacediscr} + \timestep (1 - p_2) \laplaceFactor \operator) 
		+
		\timestep^2 (\tfrac{1}{2} - p_2) 
		\br*{
			(\idmatrix{\nrspacediscr} - \timestep p_2 \laplaceFactor \operator)^{-1}
			\laplaceFactor \operator
		}^2,
	\end{split}
\end{equation}
\begin{equation}
	\label{ReactionDiffusion:eq7.1}
	\begin{split} 
		\initMatrix_{p, 2}
		=
		\timestep
		(1-\tfrac{1}{2 p_1}) 
		(\idmatrix{\nrspacediscr} - \timestep p_2 \laplaceFactor \operator)^{-1}
		+
		\timestep^2
		(\tfrac{1}{2} - p_2) 
		(\idmatrix{\nrspacediscr} - \timestep p_2 \laplaceFactor \operator)^{-1}
		\laplaceFactor \operator
		(\idmatrix{\nrspacediscr} - \timestep p_2 \laplaceFactor \operator)^{-1},
	\end{split}
\end{equation}
\begin{equation}
	\label{ReactionDiffusion:eq7.2}
	\begin{split} 
		\initMatrix_{p, 3}
		=
		\timestep
		(\tfrac{1}{2 p_1})
		(\idmatrix{\nrspacediscr} - \timestep p_2 \laplaceFactor \operator)^{-1},
		\qquad
		\initMatrix_{p, 4}
		=
		(\idmatrix{\nrspacediscr} - \timestep p_2 \laplaceFactor \operator)^{-1}
		(\idmatrix{\nrspacediscr} + \timestep (p_1 - p_2) \laplaceFactor \operator),
	\end{split}
\end{equation}
\begin{equation}
	\label{ReactionDiffusion:eq7.3}
	\begin{split}
		\initMatrix_{p, 5}
	=
		\timestep p_1 
		(\idmatrix{\nrspacediscr} - \timestep p_2 \laplaceFactor \operator)^{-1},
	\qquad
		\initMatrix_{p, 6}
	=
		\initMatrix_{p, 2} + \initMatrix_{p, 3},
	\qandq
		\initMatrix_{p, 7}
	=
		\initMatrix_{p, 5}.
	\end{split}
\end{equation}
We then assume that
for every
$p \in \params$
the parameters
$\initParams_p \in ((\R^{\nrspacediscr^d \times \nrspacediscr^d})^7)^\nrtimesteps$
are given by
\begin{equation}
\label{T_B_D}
\begin{split}
	\initParams_p = \pr{
		\underbrace{
			\initMatrix_{p}, \ldots, \initMatrix_{p}
		}_{\nrtimesteps\text{-times}}
	}.
\end{split}
\end{equation}
Roughly speaking, for every $p \in \params$ we have that 
\begin{equation}
	\label{ReactionDiffusion:eq8}
	\begin{split} 
		\adannBase_{\initParams_p}
	\approx
		\solOp
	\end{split}
\end{equation}
corresponds to an approximation of the reaction-diffusion equation in \eqref{ReactionDiffusion:eq1}
based on a finite difference discretization in space and a \LIRK\ discretization in time where the parameters $p$ correspond to the parameters of the \LIRK\ method (cf. \cref{sect:LIRK_derivation}).

\subsubsection{Numerical results for the reaction-diffusion equation}
\label{sect:RD_simul}

In this section we present numerical results for the approximation of the operator in \eqref{ReactionDiffusion:eq2}.
We test the \ADANN\ methodology with 
(see rows 13--15 in \cref{table:RD} and \cref{fig:RD_grid})
and without difference model
(see rows 10--12 in \cref{table:RD} and \cref{fig:RD_grid})
with 
	the base model and the corresponding initializations defined in \cref{sect:RD_base},
	parameter space $\params = \ReactDiffParams$,
	space discretization $\nrspacediscr = \ReactDiffSpaceStep$,
	number of time steps $\nrtimesteps \in \{2, 4, 8\}$,
	grid-based black box optimizer as described in \cref{sect:grid_optimizer},
	and
	difference model given by an \ANN\ with architecture $\ReactDiffDMarch$.
We also test different \ANN\ models with \GELU\ activation function
(see rows 1-3 in \cref{table:RD}),
\FNO\ models 
(see rows 4-6 in \cref{table:RD}),
and classical methods 
(see rows 7-9 in \cref{table:RD})
for comparison.
As classical methods we use the untrained base model
$
	\adannBase_{\initParams_{(0.5, 0.5)}}
$ 
with $\nrtimesteps \in \{2, 4, 8\}$ time steps, 
corresponding, roughly speaking, 
to a finite difference discretization in space and a CrankNicolson explicit midpoint \LIRK\ discretization in time (cf.\ \cref{sect:crank_nicolson}).
The performance of all considered methods are summarized in \cref{table:RD} and graphically illustrated in \cref{fig:RD}.
In addition, some approximations for a randomly chosen test sample are shown in \cref{fig:RD_sample}.

\begin{table} 
	\tiny
	\resizebox{\textwidth}{!}{
		\csvreader[
			tabular=|c|c|c|c|c|,
			separator=semicolon,
			table head=
			\hline 
			\thead{Method} &  
			\thead{Estimated \\ $L^2$-error} &
			\thead{Average evaluation time \\ for $\ReactDiffNrTestSamples$ test samples \\ over $\nrEvalRuns$ runs (in s)} &
			\thead{Number \\ of trainable \\ parameters} &
			\thead{Precomputation\\time (in s)}
			\\\hline,
			late after line=\\\hline
		]
		{1_numbers/rounded_methods_data_\eqname.csv}
		{
			Method=\method, 
			L2_error = \llerror, 
			nr_params = \numparams, 
			training_time = \traintime, 
			test_time = \evaltime
		}
		{\method& \llerror& \evaltime &\numparams&\traintime}%
	}
	\caption{\label{table:RD}
	Comparison of the performance of different methods for the approximation of the operator in \eqref{ReactionDiffusion:eq2} mapping source terms to terminal values of the reaction-diffusion equation in \cref{ReactionDiffusion:eq1}.}
\end{table}

\begin{figure}
	\centering
	\includegraphics[width=0.7\linewidth]{0_plots/Reaction_Diffusion/error_scatter_plot_\eqname.pdf}
	\caption{\label{fig:RD}
		Graphical illustration of the performance of the methods in \cref{table:RD}.
	}
\end{figure}

\begin{figure}
\includegraphics[width=\linewidth]{0_plots/Reaction_Diffusion/Results_grid/grid_error_overview_\eqname.pdf}
\caption{\label{fig:RD_grid}
Illustration of the \ADANN\ methodology with and without difference model (cf.\ \cref{adanns_pseudocode,adanns_pseudocode_nodiff}) applied to the approximation of the operator in \eqref{ReactionDiffusion:eq2} mapping source terms to terminal values of
the reaction-diffusion equation in \cref{ReactionDiffusion:eq1}.
\emph{Left}: Test errors of the base models prior to training as a function of the parameters used for initialization.
\emph{Middle}: Test errors of the trained base models as a function of the parameters used for initialization.
\emph{Right}: Test errors of the trained full \ADANN\ models as a function of the parameters used for initialization of the base model.
}
\end{figure}

\begin{figure}
	\minipage{0.48\textwidth}
	\expandafter\includegraphics\expandafter[width=\linewidth]{"0_plots/Reaction_Diffusion/Sample_plots/ol_plots_\eqname_0.pdf"}
	\endminipage
	\minipage{0.48\textwidth}
	\expandafter\includegraphics\expandafter[width=\linewidth]{"0_plots/Reaction_Diffusion/Sample_plots/ADANN_full_-_grid_plots_\eqname_0.pdf"}
	\endminipage
	\caption{\label{fig:RD_sample}
		Example approximation plots for a randomly chosen sample from the test set for the reaction-diffusion equation in \cref{ReactionDiffusion:eq1}.
		\emph{Left}: \ANN\ and \FNO\ approximations.
		\emph{Right}: \ADANN\ approximations.
	}
\end{figure}

\endgroup

\section{Conclusion and future work}

In this article we introduced the \ADANN\ methodology, a general framework which combines classical numerical algorithms with operator learning techniques. 
We demonstrated its effectiveness in the context of several operators related to nonlinear parabolic \PDEs, showing that the \ADANN\ methodology can outperform both classical numerical algorithms and other operator learning techniques.
There are a number of directions for further research arising from this work.
For the considered \PDE\ problems, more sophisticated and stable design algorithms for the base model could be explored.
For instance, using \LIRK\ methods of order higher than $2$ could potentially reduce the number of required time steps in the designing algorithms or allow for initial value samples with lower regularity requirements.
More generally, another natural direction is to apply the \ADANN\ methodology to other types of \PDE\ problems with different base and difference models, involving the creative design of new base models based on designing algorithms for the considered \PDE\ problems.
In that context, it would also be interesting to investigate conditions under which difference models are most effective.
Another direction is to refine the optimization over base model initializations in the \ADANN\ methodology (see \cref{sect:random_init}). This could involve developing a theoretical understanding of the objective function landscape which empirically seems to exhibit a certain regularity (cf.\ \cref{fig:SineGordon1d_grid,fig:Burgers_grid,fig:RD_grid}).
Finally, the \ADANN\ methodology seems suitable for an overall error analysis as it involves classical numerical algorithms which already have a well-established theoretical foundation.

\ifthenelse{\boolean{JMLR}}{
	\acks{This work has been partially funded by the National Science Foundation of China (NSFC) under grant number 12250610192.
	This work has also been partially funded by the Deutsche Forschungsgemeinschaft (DFG, German Research Foundation) under Germany’s Excellence Strategy EXC 2044-390685587, Mathematics M\"unster: Dynamics-Geometry-Structure.}
}{
	\section*{Acknowledgments}
	
	This work has been partially funded by the National Science Foundation of China (NSFC) under grant number 12250610192.
	This work has also been partially funded by the Deutsche Forschungsgemeinschaft (DFG, German Research Foundation) under Germany’s Excellence Strategy EXC 2044-390685587, Mathematics M\"unster: Dynamics-Geometry-Structure. 
	Moreover, this work has been supported by the Ministry of Culture and Science NRW as part of the Lamarr Fellow Network.
}

\ifthenelse{\boolean{JMLR}}{
	\appendix
}{
	\begin{appendices}
}

\crefalias{section}{appendix}
\crefalias{subsection}{appendix}

\section{Training}
\label{sect:training}

In this section we provide additional details on the training procedure used in the numerical simulations in \cref{sect:simul}.
In the training of all models we use the Adam optimizer with adaptively chosen learning rates.
Specifically, for every model and every corresponding initialization 
we approximately choose an optimal initial learning rate and,
thereafter, during the training process, successively reduce the learning rate depending on the validation error of the model.

We discuss the method by which we optimize initial learning rates in \cref{sect:learningrates} and 
we provide details on the adaptive learning rate reduction in \cref{sect:trainingdetails}.
Furthermore, we list our specific choices of hyperparameters for the experiments in \cref{sect:simul}, including the parameters governing the choice of adaptive learning rates and those for the generation of training, validation, and test sets, in \cref{table:training_hyperparams}.

\subsection{Optimal choice of initial learning rates}
\label{sect:learningrates}

In all our numerical simulations in \cref{sect:simul}, 
in the training of every model we use for every considered initialization an approximately optimal initial learning rate for the Adam optimizer.
To select such an approximately optimal initial learning rate for a model and a corresponding initialization,
we apply a golden section search (cf., e.g., \cite[Section 4.4]{Antoniou2021})
to minimize the function mapping a learning rate to
the validation error of the model
	after a fixed number of Adam training steps	with that learning rate.
The approximate minimum point identified through the golden section search is then used as the initial learning rate for the Adam training process of the considered model and initialization.
The use of the golden section search is justified by the empirical observation that the function mapping a learning rate to the validation error after a fixed number of training steps empirically tends to be unimodal, that is, there is a unique minimum and the function is strictly decreasing on the left and strictly increasing on the right of the minimum.

Choosing the initial learning rate adaptively in our experiments was motivated by the empirical observation that the standard initial learning rate of $0.001$ for the Adam optimizer often led to no improvement or even deterioration of the validation error in the training of \ADANN\ base models with highly specialized initializations.
As an example for this observation, we provide in \cref{fig:initial_LR} approximate optimal learning rates for the base models
introduced 
in \cref{SG:eq5}--\cref{SG:eq6}
in \cref{sect:SG_base}
with the initializations in \cref{SG:eq7.4} for $p=(1/2, 1/2)$.
We consider
different space discretizations (corresponding to the variable $N$ in \cref{sect:SG_base})
and different time discretizations (corresponding to the variable $M$ in \cref{sect:SG_base})
and evaluate the validation error after $50$ training steps.

We observe that the approximate optimal learning rates seem to decrease as a function of the number of space steps and the number of time steps.
Moreover, we note that the four plots on the right-hand side in \cref{fig:initial_LR} suggest that the function mapping a learning rate to the validation error after a fixed number of training steps is indeed approximately unimodal. 
Also, we observe that in the cases $N=16$, $M=64$ and $N=64$, $M=64$ (corresponding to the two plots on the lower right in \cref{fig:initial_LR}) the standard learning rate of $0.001$ for the Adam optimizer seems to lead to a clear deterioration of the validation error after $50$ training steps.

\begin{figure}[htbp]
	\centering
	\begin{subfigure}[c]{0.47\textwidth}
		\includegraphics[width=\textwidth]{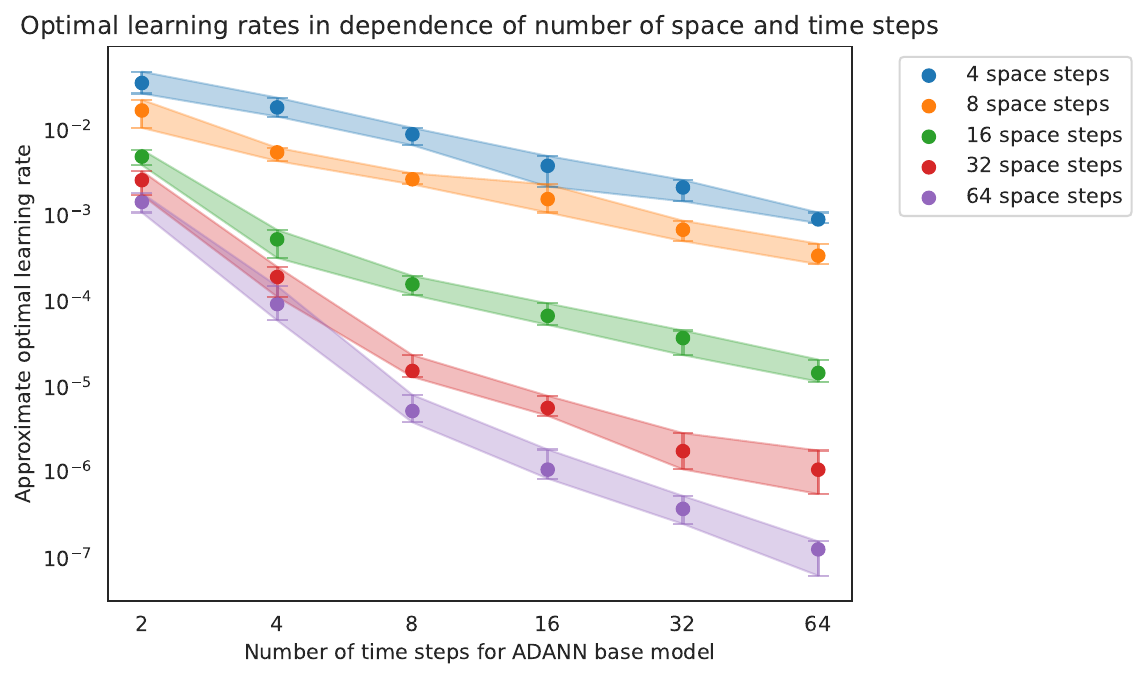}
		\label{fig:left}
	\end{subfigure}
	\hfill %
	\begin{subfigure}[c]{0.52\textwidth}
		\begin{subfigure}[c]{0.5\textwidth}
			\includegraphics[width=\textwidth]{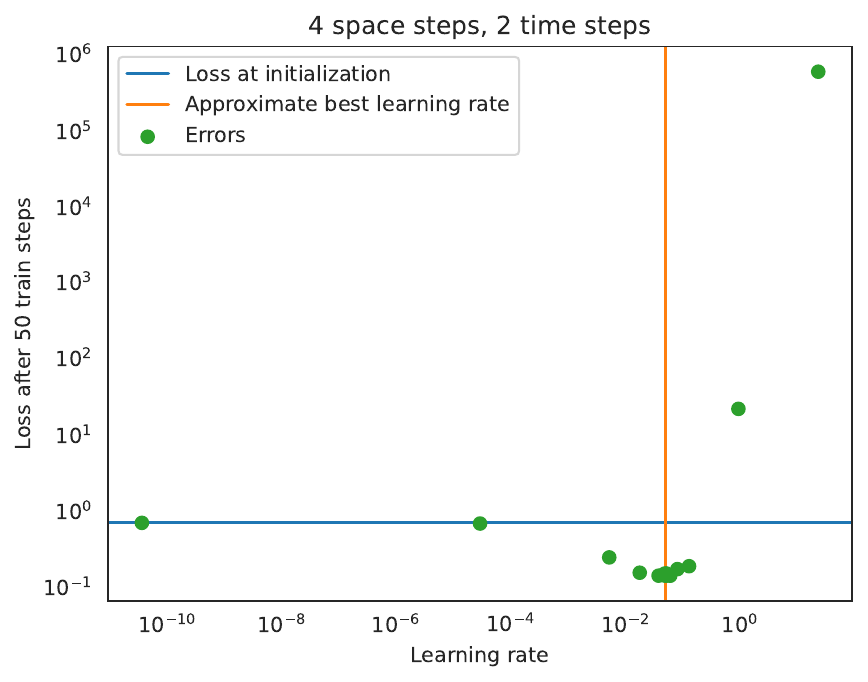}
			\label{fig:top-left}
		\end{subfigure}%
		\begin{subfigure}[c]{0.5\textwidth}
			\includegraphics[width=\textwidth]{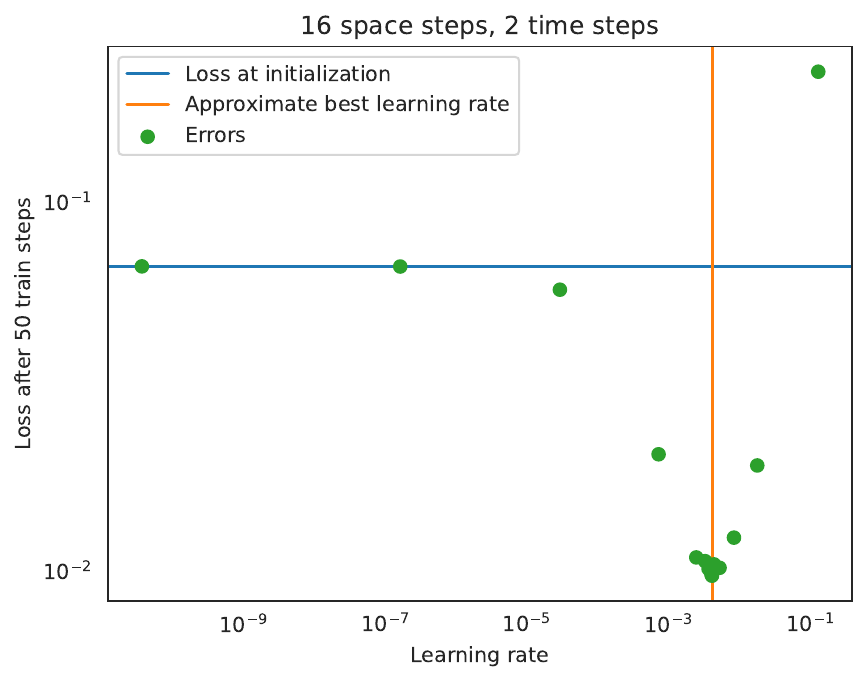}
			\label{fig:top-right}
		\end{subfigure}
		\newline %
		\begin{subfigure}[c]{0.5\textwidth}
			\includegraphics[width=\textwidth]{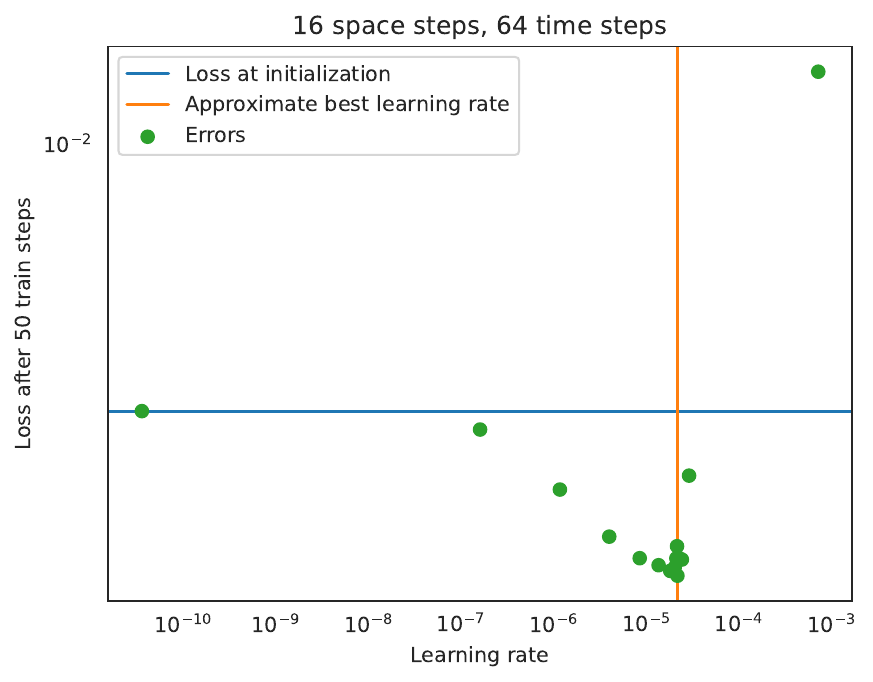}
			\label{fig:bottom-left}
		\end{subfigure}%
		\begin{subfigure}[c]{0.5\textwidth}
			\includegraphics[width=\textwidth]{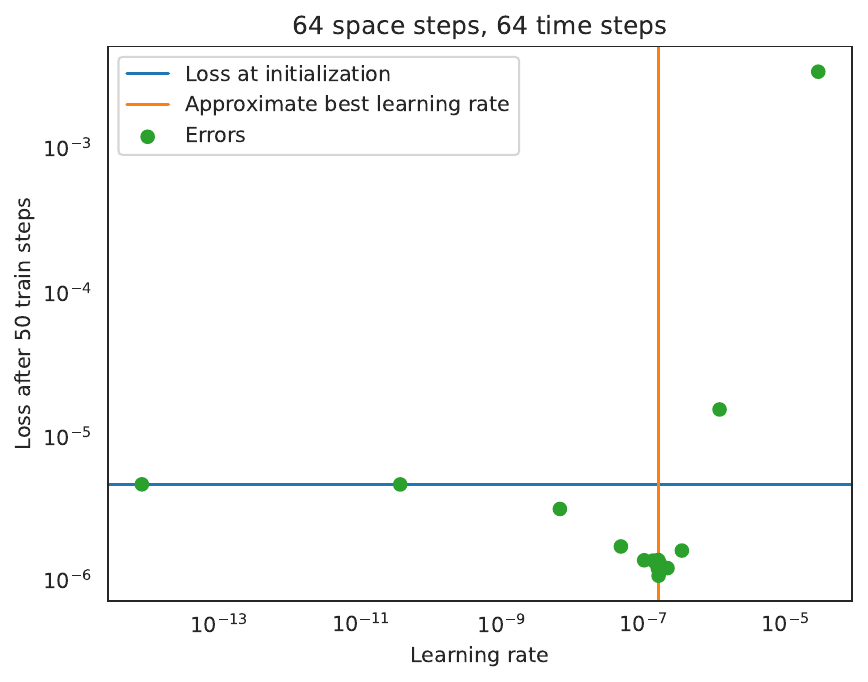}
			\label{fig:bottom-right}
		\end{subfigure}
		\label{fig:right}
	\end{subfigure}
	\caption{Optimal learning rates for the base model
	introduced in \cref{sect:SG_base}
	for the
	one-dimensional Sine-Gordon-type equation.
	\textit{Left:} Approximate optimal learning rates for different space discretizations and time discretizations in the base model. 
	The scatter points represent the average of approximated optimal learning rates from $5$ runs of golden section search
	and the error bars represent the spead over those $5$ runs.
	\textit{Right (4 plots):}
	Some examples of the validation errors after $50$ train steps computed during some of the golden section searches done for the left plot.
	}
	\label{fig:initial_LR}
\end{figure}

\subsection{Adaptive reduction of learning rates}
\label{sect:trainingdetails}

In the training of all models in the numerical simulations in \cref{sect:simul}, we adaptively reduce the learning rate and abort the training process guided by the validation error of the model being trained.
Specifically,
during the training process of a model 
we evaluate the validation error every fixed number of training steps and
whenever the relative improvement from the previous evaluation of the validation error is worse than a certain tolerance we divide the learning rate by the factor $5$.
If after a learning rate reduction, the next relative improvement of the validation error is again worse than the tolerance we abort the training process.
Within every approximation problem in \cref{sect:simul},
the number of training steps between evaluations of the validation error and the relative improvement tolerance are chosen to be the same for all considered training procedures.
The respective choices for each problem are listed in \cref{table:training_hyperparams}.

\begin{sidewaystable}[htbp]
	\centering
	\begin{adjustbox}{max width=1.1\paperwidth,center}
		\begin{tabular}{lcccc}
		\toprule
		&Sine-Gordon ($d=1$) & Sine-Gordon ($d=2$) & Burgers & Reaction-diffusion \\
		& (cf.\ \cref{sect:SineGordon1d}) & (cf.\ \cref{sect:SineGordon2d}) & (cf.\ \cref{sect:Burgers}) & (cf.\ \cref{sect:ReactionDiffusion}) \\
		\midrule
		\# space steps per dimension ($N$) 
			& $\SGonedSpaceStep$ & $\SGtwodSpaceStep$ & $\BurgersSpaceStep$ & $\ReactDiffSpaceStep$ \\
		Designing parameter set $\params$  
			& $\SGonedParams$ & $\SGtwodParams$ & $\BurgersParams$ & $\ReactDiffParams$ \\
		\midrule 
		Training &&&& \\
		\hspace{0.45cm} Batch size 
			& $\SGonedBatchSize$ & $\SGtwodBatchSize$ & $\BurgersBatchSize$ & $\ReactDiffBatchSize$ \\
		\hspace{0.45cm} \# tr.\ steps for init.\ LR search 
			& $\SGonedILRsteps$ & $\SGtwodILRsteps$ & $\BurgersILRsteps$ & $\ReactDiffILRsteps$ \\
		\hspace{0.45cm} \# tr.\ steps between val.\ error eval.\
			& $\SGonedEvalSteps$ & $\SGtwodEvalSteps$ & $\BurgersEvalSteps$ & $\ReactDiffEvalSteps$ \\
		\hspace{0.45cm} Improvement tolerance 
			& $\SGonedImprTol$ & $\SGtwodImprTol$ & $\BurgersImprTol$ & $\ReactDiffImprTol$ \\
		\midrule 
		\# of tr.\ runs per model &&&& \\
		\hspace{0.45cm} \ADANNs\ - grid
			& $\SGonedNrGridRuns$ & $\SGtwodNrGridRuns$ & $\BurgersNrGridRuns$ & $\ReactDiffNrGridRuns$ \\
		\hspace{0.45cm} \ADANNs\ - EE 
			& $\SGonedNrOptRuns$ & $\SGtwodNrOptRuns$ & $\BurgersNrOptRuns$ & $\ReactDiffNrOptRuns$ \\
		\hspace{0.45cm} \ANNs\ and \FNOs\ 
			& $\SGonedNrOLRuns$ & $\SGtwodNrOLRuns$ & $\BurgersNrOLRuns$ & $\ReactDiffNrOLRuns$ \\
		\midrule
		Algorithm for reference sol.
			& \SGonedRefAlgOne  & \SGtwodRefAlgOne  & \BurgersRefAlgOne  & \ReactDiffRefAlgOne  \\
			& \SGonedRefAlgTwo  & \SGtwodRefAlgTwo  & \BurgersRefAlgTwo  & \ReactDiffRefAlgTwo  \\
		\midrule
		Training set  &&&& \\
		\hspace{0.45cm} \# samples 
			& $\SGonedNrTrainSamples$ & $\SGtwodNrTrainSamples$ & $\BurgersNrTrainSamples$ & $\ReactDiffNrTrainSamples$ \\
		\hspace{0.45cm} \# space steps  per dimension
			& $\SGonedNrTrainSpaceDiscr$ & $\SGtwodNrTrainSpaceDiscr$ & $\BurgersNrTrainSpaceDiscr$ & $\ReactDiffNrTrainSpaceDiscr$ \\
		\hspace{0.45cm} \# time steps 
			& $\SGonedNrTrainTimeSteps$ & $\SGtwodNrTrainTimeSteps$ & $\BurgersNrTrainTimeSteps$ & $\ReactDiffNrTrainTimeSteps$ \\
		\midrule
		Validation set  &&&& \\
		\hspace{0.45cm} \# samples 
			& $\SGonedNrValidSamples$ & $\SGtwodNrValidSamples$ & $\BurgersNrValidSamples$ & $\ReactDiffNrValidSamples$ \\
		\hspace{0.45cm} \# space steps  per dimension
			& $\SGonedNrValidSpaceDiscr$ & $\SGtwodNrValidSpaceDiscr$ & $\BurgersNrValidSpaceDiscr$ & $\ReactDiffNrValidSpaceDiscr$ \\
		\hspace{0.45cm} \# time steps 
			& $\SGonedNrValidTimeSteps$ & $\SGtwodNrValidTimeSteps$ & $\BurgersNrValidTimeSteps$ & $\ReactDiffNrValidTimeSteps$ \\
		\midrule
		Test set  &&&& \\
		\hspace{0.45cm} \# samples 
			& $\SGonedNrTestSamples$ & $\SGtwodNrTestSamples$ & $\BurgersNrTestSamples$ & $\ReactDiffNrTestSamples$ \\
		\hspace{0.45cm} \# space steps per dimension
			& $\SGonedNrTestSpaceDiscr$ & $\SGtwodNrTestSpaceDiscr$ & $\BurgersNrTestSpaceDiscr$ & $\ReactDiffNrTestSpaceDiscr$ \\
		\hspace{0.45cm} \# time steps
			& $\SGonedNrTestTimeSteps$ & $\SGtwodNrTestTimeSteps$ & $\BurgersNrTestTimeSteps$ & $\ReactDiffNrTestTimeSteps$ \\
		\bottomrule
		\end{tabular}
	\end{adjustbox}
	\caption{\label{table:training_hyperparams}
	Hyperparameters for the training of the models in the numerical simulations in \cref{sect:simul}.
	}
\end{sidewaystable}

\section{Second order LIRK methods}
\label{sect:LIRK_derivation}

In this section we present a formal derivation of a well-known family of second order \LIRK\ methods for semilinear \ODEs\ (cf., e.g., \cite[Section 6.4]{Deuflhard2} and \cite{HochbruckOstermann05}) which are used to construct base models and corresponding initializations in \cref{sect:semilinear_heat,sect:simul}.
We will work in the following setting.
Let
$d \in \N$,
$\operator \in \R^{d \times d}$,
$\nonlin \in C( \R^d , \R^d)$
and consider the \ODE\
\begin{equation}
	\label{LIRK_derivation:eq1}
	\begin{split} 
		\dot{u}(t)
		=
		\operator u(t)
		+
		\nonlin(u(t))
	\end{split}
\end{equation}
for $t \in (0,\infty)$.

\subsection{Order conditions for general LIRK methods}
\newcommand{\RKsteps}{s}

We first introduce the one-step increment function of general \LIRK\ methods for the \ODE\ in \eqref{LIRK_derivation:eq1}.
Specifically, let 
$\RKsteps \in \N$,
$\alpha = (\alpha_{i, j})_{(i,j) \in \{1, 2, \ldots, \RKsteps\}^2} \in \R^{\RKsteps \times \RKsteps}$,
$\beta = (\beta_{i, j})_{(i,j) \in \{1, 2, \ldots, \RKsteps\}^2} \in \R^{\RKsteps \times \RKsteps}$,
$b = (b_i)_{i \in \{1, 2, \ldots, \RKsteps\}} \in \R^\RKsteps$
and let
$\Phi = (\Phi^h(u))_{(h, u) \in [0,\maxstep] \times \R^d} \colon [0,\maxstep] \times \R^d \to \R^d$
satisfy for all
$h \in [0,\infty)$,
$U, k_1, k_2, \ldots, k_\RKsteps \in \R^d$
with
\begin{equation}
	\label{LIRK_derivation:eq2}
	\begin{split} 
		\textstyle
		\forall \, i \in \{1, 2, \ldots, \RKsteps\}
		\colon \quad
		k_i
		=
		\operator
		(U + h \sum_{j = 1}^{i} \beta_{i, j}k_j)
		+
		\nonlin(U + h \sum_{j = 1}^{i-1} \alpha_{i, j}k_j)
	\end{split}
\end{equation}
that
\begin{equation}
	\label{LIRK_derivation:eq3}
	\begin{split} 
		\textstyle
		\Phi^h(U)
		=
		U
		+
		h
		\sum_{i = 1}^\RKsteps
		b_j k_j.
	\end{split}
\end{equation}
We refer to the number $s$ as the number of stages of the \LIRK\ method, 
we refer to $\alpha$ as the nonlinear \LIRK\ parameters,
we refer to $\beta$ as the linear \LIRK\ parameters, 
we refer to $b$ as the \LIRK\ integration weights, 
and we refer to $k_1, k_2, \ldots, k_\RKsteps$ as the \LIRK\ stages.
Although the \LIRK\ stages are defined implicitly in \eqref{LIRK_derivation:eq2}, under suitable conditions they can be computed explicitly. Specifically, under suitable conditions, we have
for all
$h \in [0,\infty)$,
$U, k_1, k_2, \ldots, k_\RKsteps \in \R^d$
with
\begin{equation}
	\label{LIRK_derivation:eq3.1}
	\begin{split} 
		\textstyle
		\forall \, i \in \{1, 2, \ldots, \RKsteps\}
		\colon \quad
		k_i
		=
		(I_d - h \beta_{i, i}\operator)^{-1}
		\pr*{
			\operator
			(U + h \sum_{j = 1}^{i-1} \beta_{i, j}k_j)
			+
			\nonlin(U + h \sum_{j = 1}^{i-1} \alpha_{i, j}k_j)
		}
	\end{split}
\end{equation}
that
\begin{equation}
	\label{T_B_D}
	\begin{split} 
		\textstyle
		\Phi^h(U)
		=
		U
		+
		h
		\sum_{i = 1}^\RKsteps
		b_j k_j.
	\end{split}
\end{equation}

Order conditions for the one-step method $\Phi$ are obtained by formally setting the Taylor expansion of $\Phi$ equal to the Taylor expansion of the solution of the \ODE\ in \eqref{LIRK_derivation:eq1} for a fixed initial value $U \in \R^d$ up to terms of a certain order.
The resulting order conditions for a second order scheme are given by
\begin{equation}
	\label{LIRK_derivation:eq4}
	\begin{split} 
		\sum_{i = 1}^\RKsteps
		b_i
		=
		1
		\qandq
		\sum_{i = 1}^\RKsteps
		b_i
		C_i
		=
		\sum_{i = 1}^\RKsteps
		b_i
		c_i
		=
		\frac{1}{2},
	\end{split}
\end{equation}
where $(C_i)_{i \in \{1, 2, \ldots, \RKsteps\}}, (c_i)_{i \in \{1, 2, \ldots, \RKsteps\}} \subseteq \R$ satisfy for all
$i \in \{1, 2, \ldots, \RKsteps\}$
that
\begin{equation}
	\label{T_B_D}
	\begin{split} 
		\textstyle
		C_i
		=
		\sum_{j = 1}^{i}
		\beta_{i, j}
		\qandq
		c_i
		=
		\sum_{j = 1}^{i-1}
		\alpha_{i, j}.
	\end{split}
\end{equation}
Under suitable regularity on the nonlinearity $\nonlin$ the conditions in \eqref{LIRK_derivation:eq4} ensure that the \ODE\ integration scheme defined through the one-step increment function $\Phi$ will have global convergence order 2.

\subsection{A family of 2 stage LIRK methods of order 2}
\label{sect:2LIRK_family}

In this section we solve the order conditions in \eqref{LIRK_derivation:eq4} in the case of $s = 2$ stages and under the assumption that $\beta_{1, 1} = \beta_{2, 2}$.
For this, let $p_1, p_2 \in (0,\infty)$ and assume that
\begin{equation}
	\label{T_B_D}
	\begin{split} 
		\alpha_{1, 2} = p_1
		\qandq
		\beta_{1, 1} = \beta_{2, 2} = p_2.
	\end{split}
\end{equation}
This and \eqref{LIRK_derivation:eq4} imply that
\begin{equation}
	\label{T_B_D}
	\begin{split} 
		b_1
		=
		1 - \tfrac{1}{2 p_1},
		\qquad
		b_2 
		=
		\tfrac{1}{2 p_1},
		\qandq
		\beta_{1, 2} = 2 p_1 (\tfrac{1}{2} - p_2).
	\end{split}
\end{equation}
This and \eqref{LIRK_derivation:eq3.1} in turn imply, under suitable conditions, that for all
$h \in [0,\infty)$,
$U, k_1, k_2 \in \R^d$
with
\begin{equation}
	\label{LIRK_derivation:eq6}
	\begin{split} 
		\textstyle
		k_1
		=
		(I_d - h p_2 \operator)^{-1}
		\pr*{
			\operator
			U
			+
			\nonlin(U)
		}
		\qand
	\end{split}
\end{equation}
\begin{equation}
	\label{LIRK_derivation:eq7}
	\begin{split} 
		k_2
		=
		(I_d - h p_2 \operator)^{-1}
		\pr*{
			\operator
			(U + h 2p_1(\tfrac{1}{2} - p_2) k_1)
			+
			\nonlin(U + h p_1 k_1)
		}
	\end{split}
\end{equation}
it holds that
\begin{equation}
	\label{LIRK_derivation:eq8}
	\begin{split} 
		\textstyle
		\Phi^h(U)
		=
		U
		+
		h
		\br[\big]{
			(1 - \tfrac{1}{2 p_1})
			k_1
			+
			(\tfrac{1}{2 p_1})
			k_2
		}.
	\end{split}
\end{equation}
We have thus derived a family of \LIRK\ methods of order two, which is parametrized by two parameters $p_1$ and $p_2$. We use this family in \cref{sect:param_LIRK}.

\subsection{The special case of the Crank-Nicolson explicit midpoint method}
\label{sect:crank_nicolson}

The scheme in \eqref{LIRK_derivation:eq6}--\eqref{LIRK_derivation:eq8} includes as a special case the well-known Crank-Nicolson explicit midpoint \LIRK\ scheme.
Specifically, note that in the special case where $p_1 = p_2 = \frac{1}{2}$ we have for all
$h \in [0,\infty)$,
$U \in \R^d$
that
\begin{equation}
	\label{T_B_D}
	\begin{split} 
		\textstyle
		\Phi^h(U)
		&=
		(I_d - \tfrac{h}{2} \operator)^{-1}
		\pr*{
			(I_d + \tfrac{h}{2} \operator)
			U
			+
			h
			\nonlin\pr*{
				(I_d - \tfrac{h}{2} \operator)^{-1}
				(U + 	\tfrac{h}{2}\nonlin(u))
			}
		}.
	\end{split}
\end{equation}

\ifthenelse{\boolean{JMLR}}{
	
}{
	\end{appendices}
}

\ifthenelse{\boolean{JMLR}}{
	\bibliography{../1_Main_bibfile/Main_bibfile.bib}
}{
	\bibliographystyle{acm}
	\bibliography{../1_Main_bibfile/Main_bibfile.bib}
}

\end{document}